# THE CONTINUOUS FUNCTIONAL CALCULUS IN LEAN

ANATOLE DEDECKER AND JIREH LOREAUX

Abstract. The continuous functional calculus is perhaps the most fundamental construction in the theory of operator algebras, especially $C^*$-algebras. Here we document our formalization of the continuous functional calculus in Lean, which constitutes the first such formalization in any proof assistant. Our implementation is already merged into Lean's mathematical library, Mathlib. We provide a brief introduction to the mathematical theory for those unfamiliar with the subject, and then highlight the design decisions in our formalization which proved to be important for usability. Our exposition is aimed at a general mathematical audience and provides a glimpse into the world of formalization by laying bare the discovery process.

## 1. Introduction

Every mathematical discipline has a toolbox which is invaluable to its practitioners, is used to lay the foundations for the subject and becomes such an essential technique that it is often taken for granted. In the theory of $C^*$-algebras, the continuous functional calculus is the first such tool. This functional calculus identifies certain unital commutative $C^*$-subalgebras with algebras of continuous functions on a compact Hausdorff space, thereby unlocking the ability to reason about $C^*$-algebras using elementary arguments about continuous functions. Our goal in this paper is to document our formalization of the continuous functional calculus in Lean in such a way that the reader can understand the motivating factors leading to this design among the multitude of possibilities. Whenever we refer to Lean herein, we always mean Lean 4.

This paper is intended to be accessible to a relatively wide mathematical audience. In particular, we hope it is accessible both to $C^*$-algebraists with little background in formalization, and also to those familiar with formalization in Lean with no background in operator algebras. We begin in Section 2 with a brief overview of the mathematical theory of the continuous functional calculus targeted at a general mathematical audience, and then discuss in Section 3 the design considerations and litmus tests which guided our formalization.

In Section 4 we walk leisurely through the design space of the continuous functional calculus in the most elementary setting — complex-valued continuous functions on the spectrum of a normal element in a unital $C^*$-algebra — accentuating the importance of usability. This presentation, while not a precise representation of the actual development process, reflects some of the approaches we considered and evaluates their shortcomings. This is intended to provide the reader, especially a mathematician not fluent in formalization, with a sense of how and why implementation in formalization matters, the questions that need to be considered, as well as how these may be satisfactorily resolved.

Section 5 piggybacks on the previous section to introduce the class-based design of the continuous functional calculus (still for *unital* $C^*$-algebras), presented as a sequence of attempts converging iteratively to the version implemented in Mathlib[1] [mC20]. While the focus of Section 4 is primarily usability, Section 5 is concerned instead with generality and flexibility.

In Section 6 we discuss the uniqueness of the continuous functional calculus, and why it is a separate class and the benefit it provides. We highlight, in terms of a few examples, a key difference in the approach to proofs utilizing uniqueness often used on paper versus in Lean. In particular, on paper, many mathematicians make direct appeal to polynomials and the Stone–Weierstrass theorem, whereas in Lean we prefer to utilize uniqueness directly. Perhaps surprisingly to the uninitiated, we develop a continuous functional calculus over $\mathbb{R}_{\geq 0}$ (the semiring of nonnegative real numbers). This has many benefits, but also contributes some unavoidable mathematical difficulties. We overcome some of them with our implementation of the uniqueness class for the continuous functional calculus, but we highlight some other difficulties in Section 11.

[1]A copy of this implementation is included in the code artifact associated to this paper.

As $C^*$-algebraists will know, non-unital algebras are an absolutely essential part of the theory. As such, our formalization would be remiss if it did not include these. While the implementation largely follows the same pattern as for unital algebras, there are a few nontrivial considerations addressed in Section 7.

Section 8 provides a brief description of the techniques we use to provide instances of all the continuous functional calculus variations we have implemented. This is somewhat interesting in its own right, because we construct instances of the non-unital calculus for $C^*$-algebras by means of the existing instances of the unital calculus, rather than a direct appeal to the Gelfand transform.

Automation is a key piece of our interface, and it is crucial to minimizing distraction to, and maximizing the productivity of, the end user. We describe this automation in Section 9. In the code artifact associated to this paper, we give example implementations of the tasks proffered in Section 3 as litmus tests to evaluate the usability of our design, thereby providing examples of our interface in action and highlighting the automation.

In Section 10 we describe a variation of our continuous functional calculus class encoding that the homomorphism is isometric, which is intentionally omitted from earlier versions. Section 11 discusses the current limitations of our approach, while Section 12 provides examples of our interface in action with example comparisons to other approaches. We conclude with Section 13 and Section 14 which list several applications of our interface that have already landed in Mathlib and some ideas for future work.

## 2. Mathematical background

A $C^*$-algebra $\mathcal{A}$ is a complex Banach algebra equipped with an involutive operation $x \mapsto x^*$, called the *adjoint*, satisfying several properties. Namely, the adjoint is *conjugate-linear* (i.e., $(a+b)^* = a^* + b^*$ and $(za)^* = \overline{z}a^*$ for any $z \in \mathbb{C}$ and $a, b \in \mathcal{A}$), *antimultiplicative* (i.e., $(ab)^* = b^*a^*$ for any $a, b \in \mathcal{A}$), and satisfies the $C^*$-identity: $\|a^*a\| = \|a\|^2$ for any $a \in \mathcal{A}$. A canonical example of a $C^*$-algebra is the algebra of $n \times n$ matrices over $\mathbb{C}$ equipped with the operator norm induced by the standard Euclidean norm, and adjoint operation given by the Hermitian adjoint. More generally, one may consider the algebra $\mathcal{B}(\mathcal{H})$ of bounded linear operators on a complex Hilbert space $\mathcal{H}$, and in fact, every $C^*$-algebra is isomorphic to a norm-closed $*$-subalgebra of the bounded operators on some Hilbert space [Dav96, Theorem I.9.12]. The standard *commutative* example of a $C^*$-algebra is the algebra $C(X, \mathbb{C})$ of continuous complex-valued functions on a compact Hausdorff space $X$, equipped with the supremum norm and the involution $*$ given by pointwise conjugation. There also exists a theory of real[2] $C^*$-algebras, but for the most part it is outside the scope of this paper, although we make brief mention of it again in Section 14.

Let us now clarify our usage of the term "algebra". In the paragraph above, and in the rest of the paper, algebras are always assumed associative, but *may or may not* contain a multiplicative identity. For instance, the ideal $\mathcal{K}(\mathcal{H})$ of $\mathcal{B}(\mathcal{H})$ consisting of *compact* endomorphisms of the Hilbert space $\mathcal{H}$, and the ideal of $C(X, \mathbb{C})$ consisting of functions vanishing at some point $x \in X$, are both examples of $C^*$-algebras with no identity elements. Following common practice in the literature on operator algebras, we will say that an algebra is **unital** if it contains an identity element, and we will always mention this explicitly. To be completely unambiguous, we will otherwise talk about **non-unital** algebras when we do not *require* the existence of an identity. A non-unital algebra may still happen to contain an identity; in the rare cases where we really want to forbid this, we will say that the algebra is **not unital**, or **without identity**. Note however that Mathlib prefers to assume algebraic structures are unital by default; hence, in the code snippets, an assumption of the form `[Algebra A]` should be read as "$A$ is a unital algebra", while `[NonUnitalAlgebra A]` corresponds to "$A$ is a non-unital algebra" in our terminology. We refer to Section 4 for more details on Lean syntax.

A key concept in the theory of $C^*$-algebras is that of the *spectrum of an element $a$* in a unital algebra $\mathcal{A}$, which is the set:

$$\sigma_{\mathbb{C}}(a) := \{\lambda \in \mathbb{C} \mid \lambda 1 - a \text{ is not invertible}\}$$

This is a purely algebraic definition, which is naturally extended to an arbitrary unital algebra $\mathcal{A}$ over a commutative ring $R$[3]. For example, if $\mathcal{A}$ is the unital algebra of $n \times n$ matrices over a field $k$, the spectrum $\sigma_k(M)$ of a matrix $M$ is precisely its set of eigenvalues over $k$. While the interest for this notion of spectrum goes far beyond the specific study of $C^*$-algebras, the spectrum enjoys some very nice properties in the setting of $C^*$-algebras. One such property is *spectral permanence*, which states that, given a unital $C^*$-subalgebra $\mathcal{B}$ of a unital $C^*$-algebra $\mathcal{A}$, an element of $\mathcal{B}$ is invertible in $\mathcal{B}$ if and only if it is invertible in $\mathcal{A}$; as a consequence, the spectrum of an element $b$ of $\mathcal{B}$ is the same whether we take $\mathcal{A}$ or $\mathcal{B}$ as the ambient $C^*$-algebra.

---

[2]Those familiar with formalization might be surprised that one would begin with the formalization of the continuous functional calculus over complex $C^*$-algebras before the formalization of real $C^*$-algebras, but it is important to realize that the real theory actually depends heavily on the complex theory. In fact, the category of real $C^*$-algebras is equivalent to the category of complex $C^*$-algebras equipped with an extra conjugate-linear *multiplicative* involution.

Sometimes the theory of $C^*$-algebras is described as noncommutative topology. The rationale behind this terminology is that *Gelfand duality* [Gel41, Satz 10, Folgergung 3], [Neg71, Pg. 228, ¶2; Pg. 237, ¶3, Theorems 2.5, 2.6, 2.8] establishes a duality (i.e., contravariant equivalence) between the category of commutative unital $C^*$-algebras and the category of compact Hausdorff spaces. More concretely, starting from a compact Hausdorff space $X$, one can build the unital $C^*$-algebra $C(X, \mathbb{C})$ of continuous complex-valued functions. This construction is contravariant in $X$, and each point of $X$ yields a *character* of this algebra, that is, a unital algebra homomomorphism from $C(X, \mathbb{C})$ to $\mathbb{C}$. This is the other functor involved in Gelfand duality: starting from a commutative unital $C^*$-algebra $\mathcal{A}$, we define its *character space* $\widehat{\mathcal{A}}$ as the set of unital algebra homomomorphisms from $A$ to $\mathbb{C}$, which are automatically continuous, of norm at most one, and $*$-preserving. Thus, $\widehat{\mathcal{A}}$ is a weak-$*$ closed subset of the unit ball in the Banach dual $\mathcal{A}^*$ of $\mathcal{A}$ so, when endowed with the weak-$*$ topology, the character space is compact Hausdorff thanks to the Banach–Alaoglu theorem. Note that this topological space may also be realized as the *space of maximal ideals* of $\mathcal{A}$, with the Zariski topology: indeed, the Gelfand–Mazur theorem asserts that maximal ideals in a unital Banach algebra are precisely the kernels of characters. Let us also mention that $\widehat{\mathcal{A}}$ is also called the *spectrum* of the algebra $\mathcal{A}$; here, the word "spectrum" connects both to the algebraic-geometry notion of spectrum (namely, a space of ideals with the Zariski topology), as well as to the notion of spectrum of an element introduced above, through the equality $\sigma_{\mathbb{C}}(a) = \{\chi(a) \mid \chi \in \widehat{A}\}$.

The last two pieces needed to build Gelfand duality are the natural isomorphisms $\Psi_X : X \to \widehat{C(X, \mathbb{C})}$ and $\Phi_{\mathcal{A}} : \mathcal{A} \to C(\widehat{\mathcal{A}}, \mathbb{C})$. The first one is given by *point evaluation*, i.e. we set $(\Psi_X(x))(f) := f(x)$ for $x \in X$ and $f \in C(X, \mathbb{C})$. To see that $\Psi_X$ is a homeomorphism, since both spaces are compact Hausdorff, it suffices to show that $\Psi_X$ is a continuous bijection. Indeed, $\Psi_X$ is continuous by the general properties of the compact-open topology, and injective by Urysohn's lemma. Surjectivity of $\Psi_X$ follows from the study of maximal ideals in $C(X, \mathbb{C})$: they are precisely kernels of the characters given by evaluation at some point $x$. Since the kernel of an arbitrary character is a maximal ideal, and characters are fully determined by their kernels, this yields the surjectivity. Finally, the isomorphism $\Phi_{\mathcal{A}} : \mathcal{A} \to C(\widehat{\mathcal{A}}, \mathbb{C})$ is the *Gelfand transform*, given by $(\Phi_{\mathcal{A}}(a))(\psi) := \psi(a)$ for $a \in \mathcal{A}$ and $\psi \in \widehat{\mathcal{A}}$. Since characters are automatically $*$-preserving, this is clearly a unital $*$-algebra homomorphism. It is norm-preserving thanks to *Gelfand's formula* for the spectral radius and the $C^*$-identity, hence it is injective and has closed range. But any two distinct characters differ on some element $a \in \mathcal{A}$, so the range of $\Phi_{\mathcal{A}}$ is a closed unital $*$-subalgebra of $C(\widehat{\mathcal{A}}, \mathbb{C})$ which separates points. By the Stone–Weierstrass theorem, it follows that the Gelfand transform is surjective.

Combining spectral permanence and Gelfand duality, we are now ready to define the *continuous functional calculus* in a unital $C^*$-algebra $\mathcal{A}$. Given an element $a \in \mathcal{A}$, we say $a$ is *normal* if it commutes with its adjoint: $a^*a = aa^*$. The unital $C^*$-subalgebra $C^*_1(a)$ generated by a normal element $a \in \mathcal{A}$ is commutative, and so by Gelfand duality above, it is isomorphic to $C(\widehat{C^*_1(a)}, \mathbb{C})$. But a character on $C^*_1(a)$ is completely determined by its values on $a$, so the space $\widehat{C^*_1(a)}$ is homeomorphic to $\{\chi(a) \mid \chi \in \widehat{A}\} = \sigma_{\mathbb{C}}(a)$, where spectral permanence ensures that $\sigma_{\mathbb{C}}(a)$ may be considered relative to the full algebra $\mathcal{A}$ instead of the subalgebra $C^*_1(a)$. One therefore obtains a (necessarily isometric) $*$-monomorphism from $C(\sigma_{\mathbb{C}}(a), \mathbb{C})$ into $\mathcal{A}$ that sends the identity function to $a$, which is referred to as the *continuous functional calculus* for $a$. In the literature, evaluation of the continuous functional calculus for a normal element $a$ at a function $f \in C(\sigma_{\mathbb{C}}(a), \mathbb{C})$ is denoted $f(a)$. If this notation seems confusing, it is worth understanding the continuous functional calculus as the *unique* (by the Stone–Weierstrass theorem) continuous $*$-homomorphism from $C(\sigma_{\mathbb{C}}(a), \mathbb{C})$ which extends the homomorphism sending a polynomial $p[z] \in \mathbb{C}[z]$ to $p(a)$.

However, there is yet more to be said. There is the *spectral mapping theorem* [KR97, Theorem 4.4.8] : if $a$ is normal and $f$ is a continuous, $\mathbb{C}$-valued function on $\sigma_{\mathbb{C}}(a)$, then $\sigma_{\mathbb{C}}(f(a)) = f(\sigma_{\mathbb{C}}(a))$, the range of $f$. Perhaps the most crucial property of the continuous functional calculus is:

---

[3]Not only that, but the reader may be surprised to see that Mathlib only requires the scalars `R` to be a commutative *semi*ring, while `A` must be a ring. This is a minor quirk that we exploit to great benefit later by allowing `R := ℝ≥0`.

**Theorem** (Composition property [KR97, Theorem 4.4.8])**.** *Let $\mathcal{A}$ be a unital $C^*$-algebra. Assume $a$ is a normal element of $\mathcal{A}$, $f$ is a continuous, $\mathbb{C}$-valued function on $\sigma_{\mathbb{C}}(a)$, and $g : \mathbb{C} \to \mathbb{C}$ is continuous on $\sigma_{\mathbb{C}}(f(a))$. Then $g \circ f$ is continuous on $\sigma_{\mathbb{C}}(a)$ and $(g \circ f)(a) = g(f(a))$.*

The notational simplicity belies its nontrivial mathematical content: on the left-hand side is the evaluation of the continuous functional calculus for $a$ at $g \circ f$, whereas on the right-hand side is the evaluation of the continuous functional calculus for $f(a)$ (which is normal) at $g$. This property depends crucially on the uniqueness of the continuous functional calculus for a given normal element $a$.

In addition, if $a \in \mathcal{A}$ is *selfadjoint* (i.e., $a^* = a$), then the spectrum of $a$ is contained in $\mathbb{R}$. For such elements, the $C^*$-algebraist often restricts attention to real-valued functions of a real variable, in which case the images of the functional calculus for $a$ are also selfadjoint. Likewise, if $a \in \mathcal{A}$ is *nonnegative* (i.e., $a = b^*b$ for some $b \in \mathcal{A}$), then the spectrum of $a$ is contained in $[0, \infty)$ (herein denoted $\mathbb{R}_{\geq 0}$), and one often restricts attention to nonnegative functions. The considerations of the previous two paragraphs allow, for instance, to define the square root of a nonnegative element $a \in \mathcal{A}$ and to conclude that $\sqrt{a^2} = a$.

The interplay between unital and non-unital $C^*$-algebras is a deep and important part of the subject, for example in the $K$-theory of $C^*$-algebras, and ensuring non-unital $C^*$-algebras are on equal footing with unital ones is an essential aspect of its formalization. Although it may seem backwards to those unfamiliar with $C^*$-algebra theory, the non-unital theory is often developed as a consequence of the unital theory. To understand why this order of development is mathematically justified, note first that one would develop the theory of compact Hausdorff spaces before proceeding to develop the theory of *pointed* compact Hausdorff spaces. Since Gelfand duality is a duality between the category of unital (respectively, non-unital) commutative $C^*$-algebras and the category of (pointed) compact Hausdorff spaces, it makes sense that one often contemplates the theory of unital $C^*$-algebras before the non-unital theory. We will return to non-unital $C^*$-algebras in Section 7.

## 3. Design considerations

In formalization, especially for constructions or theorems that experience frequent use or serve as the foundation of a wider theory, it is often the case that the seemingly most natural statement is not the one that is dictated by utilitarian design. This is quite different from the situation in much of the literature, where the focus is on establishing theorems satisfying various properties. On paper, we have the luxury of changing our perspective on the fly as we realize how to use these theorems as tools. In contrast, in formalization we must quite early on focus on the design of the interface, or else face the consequences of poor decisions as we develop the subject. In essence, the entirety of the this paper is devoted to recognizing the continuous functional calculus as a tool for future work, and designing the interface to meet the needs of that work accordingly, rather than as a theorem to be established as an end unto itself.

In this section we outline the considerations that guided the design of our interface to the continuous functional calculus in Lean, along with a few tests used to determine its suitability. For the sake of brevity, in this section we will not delve into the meaning of jargon like *unbundling* or *rewriting*. Readers unfamiliar with formalization may wish to skim this section and trust that the relevant ideas and reasoning for these design constraints will be adequately explained *in situ* in the following sections.

### 3.1. **Requirements.**

There are a number of design constraints which we placed upon ourselves in order to ensure that our interface would meet the needs of Mathlib users. Our design should:

(1) Allow for easy rewriting of terms involving the continuous functional calculus.
(2) Allow for easy use of unbundled functions (i.e., terms of type $\mathbb{C} \to \mathbb{C}$ as opposed to terms of type $C(\sigma_{\mathbb{C}}(a), \mathbb{C})$ where the continuity on the spectrum is baked into the type itself), especially in cases where those functions are obviously continuous on the spectrum; for instance, when the spectrum is finite or discrete.
(3) Provide easy and straightforward use of the composition property (i.e., proving things like $f(-a) = (x \mapsto f(-x))(a)$ and $f(a^*) = (x \mapsto f(\bar{x}))(a)$ should be nearly trivial).
(4) Ensure the full strength of the continuous functional calculus (i.e., the $*$-isomorphism with the $C^*$-subalgebra generated by the element) is recoverable from whatever implementation we choose.
(5) Allow for multiple scalar (semi)rings (namely, $\mathbb{C}$, $\mathbb{R}$ and $\mathbb{R}_{\geq 0}$) with a unified interface. This allows users to work with the continuous functional calculus on selfadjoint or nonnegative elements by considering only functions $\mathbb{R} \to \mathbb{R}$ or $\mathbb{R}_{\geq 0} \to \mathbb{R}_{\geq 0}$, respectively.

(6) Avoid requiring the algebra to be a metric space to allow for types that may not be equipped with a metric structure globally in Mathlib. The canonical use case here is matrices, which are intentionally not equipped with a metric structure in Mathlib so as to avoid choosing a preferred metric[4]. However, another comon use case is for type synonyms of $C^*$-algebras equipped with weaker topologies; for instance, the weak operator or strong operator topologies on the space of continuous linear maps, or the ultraweak topology on a von Neumann algebra.
(7) Provide a non-unital continuous functional calculus and instantiate it in both unital and not unital algebras. This allows users to define, for example, the square root of a nonnegative element in terms of the non-unital continuous functional calculus and the definition will be applicable wether the ambient algebra is unital or not.
(8) Allow for future generalization to real $C^*$-algebras.

Items (1) to (4) are all addressed to some extent in Section 4, although item (3) makes another appearance in Section 6. Items (4) to (6) are addressed in Sections 5 and 6. Item (7) is the primary focus of Sections 7 and 8. Sections 9 and 10 are related to item (6) and to items (1) and (2), respectively. Section 14 contains a brief discussion of item (8), although it is largely outside the scope of this paper.

With respect to item (5), we will introduce some terminology which is unused in the literature, but which we find useful for disambiguation. When we refer to a (possibly non-unital) continuous functional calculus *over* $R$, where $R$ is a commutative semiring, we mean a continuous functional calculus which allows for functions $f : R \to R$. This will only be applicable when the algebra $A$ is an $R$-algebra and the elements for which such a calculus is valid depends on $R$ (i.e., for $\mathbb{R}$ or $\mathbb{R}_{\geq 0}$, one must require the elements of the algebra to be selfadjoint or nonnegative, respectively).

There are several tasks which the (unital or non-unital) continuous functional calculus should be able to accomplish, and they require essentially no prerequisites. Before embarking on this project, we set for ourselves a series of litmus tests to evaluate whether our design was sufficiently usable to be effective in practice. If these were exceedingly difficult to implement, then the design was not suitable. Example implementations of each can be found in the code artifact associated to this paper, except for the last one which is nevertheless implemented in Mathlib already.

**Inverses:** For an invertible element $a$ in a unital $C^*$-algebra and a function $f$ which is continuous on the spectrum of the inverse of $a$, show that: $f(a^{-1}) = (x \mapsto f(x^{-1}))(a)$. This test requires:
(1) A unital continuous functional calculus over $\mathbb{C}$ for normal elements.
(2) Use of functions which are not everywhere continuous[5] (i.e., $x \mapsto x^{-1}$).
(3) Application of the composition property (since we are working with continuous functional calculi for both $a$ and $a^{-1}$).

**Positive and negative parts:** For a selfadjoint element $a$ in a non-unital $C^*$-algebra, there is a unique pair $(a^+, a^-)$ of commuting positive elements such that $a = a^+ - a^-$ and $a^+a^- = 0$. This test has slightly different requirements:
(1) A non-unital continuous functional calculus over $\mathbb{R}$ for selfadjoint elements.
(2) Use of automatically inserted tactics to prove $f(0) = 0$ for various functions $f$. In Lean, this insertion is accomplished by an automatic parameter, or `autoParam`. The tactic to use for a given hypothesis is specified when a lemma is declared, and when the lemma is later used, Lean calls this tactic to satisfy the hypothesis.

**Square roots:** For a nonnegative element $a$ in a non-unital $C^*$-algebra, there is a unique nonnegative element $a^{\frac{1}{2}}$ so that $(a^{\frac{1}{2}})^2 = a$. This requires:
(1) A non-unital continuous functional calculus over $\mathbb{R}_{\geq 0}$ for nonnegative elements.
(2) Application of the composition property, which requires uniqueness of the functional calculus. This is trickier over $\mathbb{R}_{\geq 0}$ than over $\mathbb{R}$ or $\mathbb{C}$ because the Stone–Weierstrass theorem doesn't apply.
(3) Showing that we can use different functional calculi (i.e., over $\mathbb{R}$ too!). To be more precise, given a nonnegative element $a$, we can define the square root of $a$ using either the functional calculus over $\mathbb{R}_{\geq 0}$ with the function `NNReal.sqrt` or the functional calculus over $\mathbb{R}$ with the function `Real.sqrt`, and we can show these are equal. This matters if, for example, if we want to consider $f(\sqrt{a})$ for some $f : \mathbb{R} \to \mathbb{R}$ and use the composition property to realize it as $(x \mapsto f(\sqrt{x}))(a)$.

[4]The reasoning here is that $n \times n$ matrices of complex numbers represent the endomorphisms of any $n$-dimensional complex normed space, or even the Hilbert-Schmidt or trace-class operators of a $n$-dimensional complex Hilbert space. Any of these choices would give rise to a different choice of a natural norm, and making a default choice would cause troubles for users wishing to work with another one.

**Matrices over $\mathbb{R}$:** Is it possible to create an instance of the continuous functional calculus on selfadjoint elements of `Matrix n n ℝ`, i.e., $n \times n$ matrices over $\mathbb{R}$? The primary reason to include this test is that it requires the continuous functional calculus to work for algebras over $\mathbb{R}$ that are not equipped with a metric structure.

## 4. An interface for unital $C^*$-algebras

The aim of this section is to motivate and describe the most fundamental design choices made in our formalization: *unbundling*, and the use of so-called *"junk values"*. In order to keep the presentation simple and, hopefully, accessible to a broad mathematical audience with no background in formalization, we focus on the most basic setting, that of normal elements in *unital* complex $C^*$-algebras, and we postpone the detailed exposition of more technical requirements to later sections. Because we restrict attention to this most elementary setting, we emphasize that the exposition in this section does not actually represent faithfully the design process, and should rather be understood as a motivation and description of some general principles.

4.1. **A mathematician's definition.** The first step of our work was of course to formalize the mathematical content underlying the continuous functional calculus, namely the (unital) Gelfand transform and spectral permanence. In fact, these were added to Mathlib back in February 2023, together with the usual definition of continuous functional calculus described in Section 2. Although we won't focus on the formalization of the prerequisites, we recall the first definition of the functional calculus in Listing 1, which will serve as our starting point for the design process. Let us also use it to explain some glimpses of Lean syntax for the unfamiliar reader.

First, the `variable` line declares `A` to be an arbitrary unital $C^*$-algebra. In that context, we **def**ine a new *term*, i.e a mathematical object, named `continuousFunctionalCalculus`, which takes as input a normal element `a` of `A`, and outputs a term of *type* `C(spectrum ℂ a, ℂ) ≃⋆ₐ[ℂ] StarAlgebra.elemental ℂ a`. Here, the expression `C(spectrum ℂ a, ℂ)` denotes the type of continuous functions on the spectrum of `a` with values in `ℂ`, while `StarAlgebra.elemental ℂ a` is a name for the the unital $C^*$-subalgebra generated by $a$, and the notation `≃⋆ₐ[ℂ]` refers to the type of $*$-isomorphisms between these $*$-algebras (over `ℂ`). It is important for mathematicians not well-versed in formalization to be aware that this constitutes a **def**inition in Lean (as opposed to a `theorem`) because this term contains *data*, namely the functions which constitute the $*$-isomorphism. This therefore contains nontrivial mathematical content, not simply terminology to be used later. The actual *definition* of the term `continuousFunctionalCalculus a`, expressed using the `gelfandStarTransform` of the $C^*$-algebra `StarAlgebra.elemental ℂ a`, comes after the symbol `:=`. Finally, we prove a `theorem`, stating that the continuous functional calculus for some normal element `a` maps (the restriction to the spectrum of `a` of) the identity function to `a`. In Lean, it is possible to temporarily omit a proof using the `sorry` keyword, which allows Lean to accept the theorem for future use while still remembering that this theorem is not yet proven; it can even be used to temporarily omit a certain step in a proof.

Listing 1: The continuous functional calculus as a $*$-isomorphism
*LeanCFC/Snippets/Baseline.lean*

```
8  variable {A : Type*} [CStarAlgebra A]
9
10 /-- The ⋆-isomorphism between the continuous functions on the spectrum of `a`
11 and the (unital) C⋆-subalgebra generated by `a`. -/
12 def continuousFunctionalCalculus (a : A) [IsStarNormal a] :
13     C(spectrum ℂ a, ℂ) ≃⋆ₐ[ℂ] StarAlgebra.elemental ℂ a :=
14   characterSpaceHomeo a |>.compStarAlgEquiv' ℂ ℂ |>.trans
15     (gelfandStarTransform (elemental ℂ a)).symm
16
```

[5] Note that we say "not everywhere continuous" rather than "not everywhere defined" because, in Mathlib, $0^{-1} = 0$ when 0 is the additive identity in a field. So indeed the function $x \mapsto x^{-1}$ is everywhere defined, but not continuous at 0. For a friendly explanation of some reasons for the $0^{-1} = 0$ choice, see [Buz20]

```
17 /-- The continuous functional calculus identifies the (restriction of) the
18 identity function on the spectrum of `a` with `a`. -/
19 theorem continuousFunctionalCalculus_map_id (a : A) [IsStarNormal a] :
20     continuousFunctionalCalculus a ((ContinuousMap.id ℂ).restrict (spectrum ℂ a)) =
21       ⟨a, self_mem ℂ a⟩ :=
22   gelfandStarTransform (elemental ℂ a) |>.symm_apply_apply _
```

While this definition is completely satisfying from a mathematical point of view, let us recall that the functional calculus is not the end of the story. Rather, it is a crucial stepping stone on the road towards the general theory of $C^*$-algebras, where it will be used extensively. Thus, to progress towards formalization of these later steps, it is essential to provide an interface dictated by utilitarian design rather than mathematical simplicity. Listing 2 presents such an interface, which we will explain, motivate and expand in the rest of this section. Although this is not yet the version that one can find in Mathlib, it is a good approximation: more precisely, in the context of unital complex $C^*$-algebras, `cfc f a` (defined *ibid*) agrees with the one in Mathlib.

Listing 2: A more convenient interface for the continuous functional calculus
*LeanCFC/Snippets/Intermediate.lean*

```
6 variable {A : Type*} [CStarAlgebra A]
7
8 open scoped Classical in
9 def cfc (f : ℂ → ℂ) (a : A) : A :=
10   if h : IsStarNormal a ∧ Continuous ((spectrum ℂ a).restrict f)
11   then
12     letI := h.1
13     continuousFunctionalCalculus a ⟨(spectrum ℂ a).restrict f, h.2⟩
14   else 0
15
16 lemma cfc_def {f : ℂ → ℂ} {a : A} [ha : IsStarNormal a] (hf : ContinuousOn f (spectrum ℂ a)) :
17     cfc f a = continuousFunctionalCalculus a ⟨_, hf.restrict⟩ := by
18   rw [continuousOn_iff_continuous_restrict] at hf
19   simp [cfc, ha, hf]
20
21 lemma cfc_id (a : A) [ha : IsStarNormal a] :
22     cfc id a = a := by
23   rw [cfc_def continuousOn_id]
24   exact congrArg _ (continuousFunctionalCalculus_map_id a)
```

4.2. **Simple expressions.** We'll now compare Listings 1 and 2 to explain and highlight some conveniences afforded by the latter. The distinction between these listings is apparently completely elementary: we merely extend[6] the definition of $f(a)$ to any map $f : \mathbb{C} \to \mathbb{C}$ and any element $a \in A$, by letting $f(a) = 0$ whenever $a$ is not normal or $f$ is not continuous on $\sigma_{\mathbb{C}}(a)$. Mathematically, this is completely uninteresting, since none of the properties — not even linearity in the variable $f$! — extend without these hypotheses. Nevertheless, this kind of extension by "junk values" is a common practice in formalization, perhaps the most prevalent example being that Mathlib, as well as most proof assistants, defines `(0 : ℝ)⁻¹ := 0`.

The motivation behind this is that we want to emulate the "Write first, think later" practice of informal mathematics. On paper, it is very common to write down an entire, complicated computation, before checking (or letting the reader check) that every term is well-defined — e.g, each denominator is nonzero, each integrand is integrable, and so forth — and that the manipulations of these terms are legal. While Lean is prefectly fine with waiting a bit for a proper justification of some proof step, it has no tolerance for ill-defined expressions: after all, the core job of the Lean kernel is precisely to check that each expression, whether it denotes a mathematical object or a proof, is well-formed. By extending the set of well-defined expressions using junk values, we replace ill-*defined* objects by ill-*behaved* ones, which means that the usual checks are delayed until we have to prove something nontrivial about these — e.g, `a⁻¹ * a = 1` only holds for `a ≠ 0`.

[6]One could argue that we also restrict our attention to functions defined on all of $\mathbb{C}$, but since we don't care about continuity outside of $\sigma_{\mathbb{C}}(a)$ this doesn't change anything mathematically.

Coming back to functional calculus, it may not be clear how this discussion applies to Listings 1 and 2. After all, the expression `((ContinuousMap.id ℂ).restrict (spectrum ℂ a))` does not seem to contain any proof of normality or continuity. In the case of continuity, this is because the expression for the identity map `((ContinuousMap.id ℂ).restrict (spectrum ℂ a))` has type `C(spectrum ℂ a, ℂ)`, meaning that it comes *bundled* [Baa22, Section 6] with a proof of continuity. Unfortunately, this only works because Mathlib provides us with the relevant building blocks — the identity as an element of `C(ℂ, ℂ)`, and a restriction operation specific to continuous functions. By contrast, the expression $\exp(a)$ would be written as in Listing 3, where `by fun_prop` is a proof of continuity. Regarding normality, Listing 1 relies on *typeclass inference*: the square brackets around `IsStarNormal a` indicate to Lean that it should search for a normality proof using the *instances* declared in the library. For example, Lean will be able to *synthesize* (i.e., find) an instance of `IsStarNormal a` when the algebra `A` is commutative. While this is practical in basic cases, one cannot always rely on the typeclass system. This happens for example in the expression[7] $f(a+a)$: since it is not true that the sum of two normal elements is normal, Lean doesn't know that $a+a$ is normal. Thus, to write down such an expression using Listing 1, one would need to first prove normality of $a+a$[8]. This is demonstrated in Listing 4.

Listing 3: Definition of $\exp(a)$ using the approach of Listing 1
*LeanCFC/Snippets/Baseline.lean*

```
24 example {a : A} [IsStarNormal a] : StarAlgebra.elemental ℂ a :=
25   continuousFunctionalCalculus a ⟨fun x ↦ Complex.exp x, by fun_prop⟩
```

Listing 4: Definition of $f(a+a)$ using the approach of Listing 1
*LeanCFC/Snippets/Baseline.lean*

```
27 example {a : A} [ha : IsStarNormal a] {f : C(spectrum ℂ (a + a), ℂ)} :
28     StarAlgebra.elemental ℂ (a + a) :=
29   have : IsStarNormal (a + a) := by
30     rw [isStarNormal_iff, star_add]
31     exact (ha.star_comm_self.add_left ha.star_comm_self).add_right
32       (ha.star_comm_self.add_left ha.star_comm_self)
33   continuousFunctionalCalculus (a + a) f
```

Another convenient feature of Listing 2 is that it avoids the use of *subtypes* like `StarAlgebra.elemental ℂ a` or `spectrum ℂ a`, which are subtypes of `A` and `ℂ`, respectively. To explain what we mean by "subtype", we need to point out a key distinction between Lean's underlying theory and a more familiar set theory. In set theory, the $\in$ relation is nothing more than a binary relation on all mathematical objects (i.e., sets). In particular, any set can (and does) belong to a lot of different sets, and there is no preferred one. By contrast, in Lean, each term has a *unique* type, hence there is no notion of one type being a subset of another. Instead, what one can do is declare canonical maps between types, called *coercions* or *casts*, which Lean uses to interpret elements of some type as elements of another type. Lean tries to insert these automatically, but one easily gets into situations where there is not enough information to do so. In these cases, the user needs to provide a hint to help Lean determine it should insert the coercion, which reduces readability. Such hints take the form of *type ascriptions* where the user explicitly indicates the expected type of the term. As an example, a user might specify `(n : ℝ)` when `n` has type `ℕ`; this causes Lean to check that `n` indeed has type `ℝ`, and upon concluding that it *doesn't*, Lean then searches for a coercion from `ℕ` to `ℝ`. Even when inserted automatically without type ascriptions from the user, Lean must keep track of these functions in the expressions, where they appear as `↑`. This means one has to constantly use compatibility results like `↑(a + b) = ↑a + ↑b` or `↑0 = 0`, which tend to make computations messy.

[7] If that looks like something no one would ever write, remember that it could just be a simple step in the middle of a complicated computation.

[8] To be complete, there are still ways to work around proving normality manually by declaring carefully-chosen instances. Another way would be to compute `a + a` in the commutative unital $C^*$-subalgebra generated by $a$. See Listing 11 and the surrounding discussion for an example.

For example, the appearance of the type `spectrum ℂ a` in Listing 1 implies that we constantly have to restrict `f` to the spectrum of `a`. This pain point may already be observed in Listing 1 in the case of the identity function, and Listing 5 shows how this gets even worse when `f` is defined on another subtype which "contains" the spectrum — keep in mind these complicated expressions are just for being able to *write* $f(a)$! Of course, the approach in Listing 2 does suffer from the dual issue, in the sense that one may need to extend a function to be defined on all of $\mathbb{C}$. Here though, the common practice of extending functions using junk values works in our favor, since it means that essentially any function one might want to write down is *already* defined everywhere.

Listing 5 concerns a unitary element $u$ (i.e., $uu^* = 1 = u^*u$) in a unital $C^*$-algebra $\mathcal{A}$. It is an elementary fact that the spectrum $\sigma_{\mathbb{C}}(u)$ is contained in the unit circle $\mathbb{S}^1$ when $u$ is unitary. Indeed, the $C^*$-identity implies that $\|u^{-1}\| = \|u^*\| = 1 = \|u\|$, and the spectral radius $\rho(u)$ is bounded by $\|u\|$, so the spectrum is contained in the closed unit disk. But because $\|u^{-1}\| = 1$ and $\sigma_{\mathbb{C}}(u^{-1}) = \sigma_{\mathbb{C}}(u)^{-1}$, the spectrum must also be contained in the complement of the open unit ball, and therefore is a subset of the unit circle. Consequently, it is common to apply functions which are continuous on the unit circle to unitary elements, as we do below. This listing illustrates that simply writing $f(u)$ for $f \in C(\mathbb{S}^1, \mathbb{C})$ is nontrivial using Listing 1. Note that the expression for $f(u)$ (which occurs after the `:=`) is actually more complicated than it appears, since Lean automatically inserts the coercion from `StarAlgebra.elemental ℂ u` to `A`. Not only that, but the type ascription `(u : A)` is necessary for Lean to insert the coercion from `unitary A` to `A`.

Listing 5: Troubles with restrictions : $f(u)$ for $f \in C(\mathbb{S}^1, \mathbb{C})$ and $u$ unitary
*LeanCFC/Snippets/Baseline.lean*

```
39 example {u : unitary A} {f : C(Metric.sphere (0 : ℂ) 1, ℂ)} : A :=
40   continuousFunctionalCalculus (u : A)
41     (f.comp (.inclusion (Unitary.spectrum_subset_circle u)))
```

Before moving on, let's impart some nuance to the first comparison between Listing 1 and Listing 2. If we compare the exact mathematical content of Listing 2 to that of Listing 1, we have lost some significant information about the functional calculus: with this definition, $f \mapsto f(a)$ is neither surjective (because we forgot its image), nor injective (because it doesn't depend on the value of $f$ outside of $\sigma_{\mathbb{C}}(a)$), nor even an algebra homomomorphism (because of junk values). More subtly, we lost the fact that it is an isometry when restricted to $C(\sigma_{\mathbb{C}}(a), \mathbb{C})$, which is of crucial importance for the rest of the theory. Listing 6 shows how this can be fixed rather easily by restating all the relevant facts as separate lemmas in the language of Listing 2 (Let us note that the corresponding proofs are quite messy, precisely because we have to use the interface of Listing 1. We omit some of them for clarity).

Nevertheless, it is not clear at this stage that any of this is worth the trouble. First, setting up a huge interface for a specific construction does come with some disadvantages in terms of maintenance cost and steepening the learning curve; for example, a new user might try to use the generic `map_add` rather than the specific `cfc_add`. Furthermore, the use of bare functions instead of bundled continuous maps or bundled $*$-algebra homomomorphisms causes some friction when using the library. The statement of `norm_cfc_eq` in Listing 6 illustrates such friction: if `f` had type `C(spectrum ℂ a, ℂ)`, the right hand side could be written naturally as `‖f‖`, but the use of bare functions forces us to write it down more explicitly as the best positive bound for all the $\|f(x)\|$, $x \in \sigma_{\mathbb{C}}(a)$. Some of these inconveniences were significant enough for us to stick to bundled continuous maps and $*$-algebra homomomorphisms in our initial designs. The authors are indebted to Frédéric Dupuis for suggesting that the unbundled approach is worth the cost. We will now see more profound arguments in favor of Listing 2.

Listing 6: Recovering some properties in the setting of Listing 2
*LeanCFC/Snippets/Intermediate.lean*

```
41 lemma cfc_add {f g : ℂ → ℂ} {a : A} [ha : IsStarNormal a]
42     (hf : ContinuousOn f (spectrum ℂ a)) (hg : ContinuousOn g (spectrum ℂ a)) :
```

```
43     cfc (f + g) a = cfc f a + cfc g a := by
44   rw [cfc_def (hf.add hg), cfc_def hf, cfc_def hg]
45   norm_cast
46   exact map_add (continuousFunctionalCalculus a) ⟨_, hf.restrict⟩ ⟨_, hg.restrict⟩
47
48 lemma cfc_neg {f : ℂ → ℂ} {a : A} [ha : IsStarNormal a]
49     (hf : ContinuousOn f (spectrum ℂ a)) :
50     cfc (-f) a = -(cfc f a) := by
51   rw [Pi.neg_def, cfc_def hf, cfc_def hf.neg]
52   norm_cast
53   exact map_neg (continuousFunctionalCalculus a) ⟨_, hf.restrict⟩
54
55 -- Similarly, we prove `cfc_mul`, `cfc_sub`, `cfc_const`...
56
57 lemma cfc_eq_iff {f g : ℂ → ℂ} {a : A} [ha : IsStarNormal a]
58     (hf : ContinuousOn f (spectrum ℂ a))
59     (hg : ContinuousOn g (spectrum ℂ a)) :
60     cfc f a = cfc g a ↔ ∀ x ∈ spectrum ℂ a, f x = g x := by
61   simp [cfc_def hf, cfc_def hg, ContinuousMap.ext_iff]
62
63 lemma range_cfc_eq {a : A} [ha : IsStarNormal a] :
64     Set.range (fun f ↦ cfc f a) = StarAlgebra.elemental ℂ a :=
65   sorry
66
67 lemma norm_cfc_eq {f : ℂ → ℂ} {a : A} [ha : IsStarNormal a]
68     (hf : ContinuousOn f (spectrum ℂ a)) :
69     ‖cfc f a‖ = sInf {C : ℝ | 0 ≤ C ∧ ∀ x ∈ spectrum ℂ a, ‖f x‖ ≤ C} :=
70   sorry
```

4.3. **Rewriting with ease by avoiding dependent types.** Following computer science terminology, *rewriting* refers to the act of replacing one expression `e` by an *equal* expression `e'`, inside some larger statement or expression. In Lean this is accomplished by using the tactic `rw` (short for "rewrite"): if `H` is a proof of some equality `a = b`, `rw [H]` replaces all occurences of `a` by `b` in the goal or in some other specified location. This tactic underlies essentially every computation in Lean, and as such one should expect to be very frequently rewriting either $f$ or $a$ in an application $f(a)$ of the functional calculus.

During formalization, utilizing this substitution property of equality employed by rewriting is effiicient, and close to standard mathematical practice. Moreover, due to the fundamental nature of equality in Lean, rewriting works out of the box in most settings. To understand where things can go wrong, we need to talk very briefly about the way rewriting works under the hood, although a precise explanation is far beyond the scope of this article. Essentially, when trying to rewrite `a = b`, with `a` and `b` of type `α`, inside an expression `e`, Lean searches for some function `φ`, called the *motive*, such that `e` coincides with `φ a`. Rewriting `a = b` then simply amounts to the fact that a function maps equal arguments to equal values. The potential issue here is that the function `φ` has to be defined, and well-typed, for *all* elements of type `α`, not just `a` and `b`. When this is not the case, the user will encounter error messages of the form "Motive is not type correct…".

This observation is enough to suggest that rewriting `a = b` inside `continuousFunctionalCalculus a f` will rarely work! Indeed, there are *two* subterms in this expression whose type depend on `a`: the implicit parameter of type `IsStarNormal a` (found by typeclass inferance), and the function `f` of type `C(spectrum ℂ a, ℂ)`. Thus, when Lean tries to find a motive `φ` mapping `a` to `continuousFunctionalCalculus a f`, it is faced with two issues: the proof of normality for `a` may not be generalizable to a proof of normality for an arbitrary element `x` (and in fact, such an element may fail to be normal), and there may be no way to reinterpret `f` as an element of `C(spectrum ℂ x, ℂ)`.

In fact, we can illustrate these two points of failure separately. Firstly, in Listing 7, we apply the functional calculus to functions of the form `f.restrict (spectrum ℂ _)`. In other words, these are *syntactically* restrictions of some globally defined continuous function `f` to the spectrum of some element. Thus, `rw` has no trouble generalizing `a` to an arbitrary element `x` in this part of the expression. The proofs of normality, however, are specific to `a` and `b`, hence there is no way to find some `h : ∀ x, IsStarNormal x` abstracting `ha` and `hb`. This prevents `rw` from finding a type-correct motive.

Listing 7: Rewriting failures in the setting of Listing 1
*LeanCFC/Snippets/Baseline.lean*

```
70 example {a b : A} [ha : IsStarNormal a] [hb : IsStarNormal b] (hab : a = b) {f : C(ℂ, ℂ)}
71     (H : continuousFunctionalCalculus a (f.restrict (spectrum ℂ a)) = 0) :
72     continuousFunctionalCalculus b (f.restrict (spectrum ℂ b)) = 0 := by
73   rw [hab] at H -- "motive is not type correct" because of `IsStarNormal a`
```

Conversely, in Listing 8, we assume that the ambient unital $C^*$-algebra is *commutative*. This guarantees not only that all elements are normal, but that they are all normal "for the same reason". Concretely, all the proofs of normality in this statement are given by some general lemma applied to an element, hence Lean can figure out what the proof should look like for an arbitrary element. Note that it is really the *form* of the proof which matters here: if we add `[ha : IsStarNormal a]` then Lean fails again to infer this part of the motive. To illustrate the second point of failure, this time we apply the functional calculus to a function `f : C(spectrum ℂ a, ℂ)`, which cannot be interpreted as having type `f : C(spectrum ℂ x, ℂ)` for an arbitrary element `x`. We can still make sense of the statement, by using `cast (by rw [hab]) f` (which is essentially "`f` reinterpreted as having type `C(spectrum ℂ b, ℂ)`"), but this is cumbersome and does not help `rw` to build a motive. More specifically, the issue is now that the *proof of equality* `hab` cannot be generalized to an arbitrary element.

Listing 8: Rewriting failures in the setting of Listing 1
*LeanCFC/Snippets/Baseline.lean*

```
103 example {B : Type*} [CommCStarAlgebra B] {a b : B} (hab : a = b) {f : C(spectrum ℂ a, ℂ)}
104     (H : continuousFunctionalCalculus a f = 0) :
105     continuousFunctionalCalculus b (cast (by rw [hab]) f) = 0 := by
106   rw [hab] at H -- "motive is not type correct" because of `spectrum ℂ a`
```

More generally, these issues tend to appear when a function takes arguments whose type depends on a former argument. As the approach in Listing 2 avoids types depending on `a` such as `spectrum ℂ a`, `IsStarNormal a`, and `StarAlgebra.elemental ℂ a`, rewriting works out of the box in the expected way, as demonstrated in Listing 9.

Listing 9: Rewriting in the setting of Listing 2
*LeanCFC/Snippets/Intermediate.lean*

```
72 example {a b : A} [IsStarNormal a] [IsStarNormal b] (hab : a = b) {f : ℂ → ℂ}
73     (H : cfc f a = 0) : cfc f b = 0 := by
74   rw [hab] at H
75   exact H
```

4.4. **The composition property and its consequences.** The issues mentioned in the previous section already strongly suggest to abandon the approach of Listing 1. However, in our development process, it was actually consideration of the composition property which clinched the switch to Listing 2. Together with the fact that $f \mapsto f(a)$ is an isometric $*$-morphism, the composition property is probably the most crucial property for applications of the functional calculus, as well as the only way to link the functional calculi associated to distinct elements of the algebra. Thus, convenient use of this property and its consequences is of primary importance.

In this regard, the approach of Listing 1 fails quite badly, as even *stating* the composition property is convoluted. Indeed, for the expression $g(f(a))$ to make sense, we need $f \in C(\sigma_{\mathbb{C}}(a), \mathbb{C})$ and $g \in C(\sigma_{\mathbb{C}}(f(a)), \mathbb{C}) = C(f(\sigma_{\mathbb{C}}(a)), \mathbb{C})$, and for such $f$ and $g$ one cannot formally write $g \circ f$. Instead, one has to consider the *co-restriction* of $f$ as a function $\hat{f} : C(\sigma_{\mathbb{C}}(a), \sigma_{\mathbb{C}}(f(a)))$, to be able to write $g(f(a)) = (g \circ \hat{f})(a)$. Notice in particular that this adds even more types depending on $f$ and $a$, which means that the issues raised in the previous section could get even worse. On top of that, $\hat{f}$ is only well-defined because of the spectral mapping theorem, which would thus need to be invoked repeatedly, and the co-restriction adds complexity to the expressions we are working with. This only gets worse with more iterations of composition, where one may need to use facts such as $h \circ \hat{g} \circ \hat{f} = h \circ \widehat{g \circ \hat{f}}$. One can get a slightly cleaner statement by asking the user to provide the co-restricted avatar of $\hat{f}$, as demonstrated in Listing 10. This does give more flexibility for rewriting since one can delay the proof that the given $\hat{f}$ indeed coincides with $f$, but this is still far from satisfactory. We encourage the reader, without accessing Lean's interface, to try to parse the statement of the composition property in Listing 10; in particular, try to consider carefully in which type the equality occurs, in which type the continuous functional calculus is applied, and why the relevant elements are normal.[9]

Listing 10: Statement of the composition property for the approach of Listing 1
*LeanCFC/Snippets/Baseline.lean*

```
109 theorem continuousFunctionalCalculus_map_comp {a : A} [IsStarNormal a]
110     {f : C(spectrum ℂ a, ℂ)}
111     {f' : C(spectrum ℂ a, spectrum ℂ (continuousFunctionalCalculus a f))}
112     {g : C(spectrum ℂ (continuousFunctionalCalculus a f), ℂ)}
113     (H : ∀ x, f x = f' x) :
114     continuousFunctionalCalculus a (g.comp f') =
115     continuousFunctionalCalculus (continuousFunctionalCalculus a f) g := by
116   sorry
```

By contrast, stating and using the composition property in the approach of Listing 2 is as simple as it could be: since all functions are defined on the whole complex plane, one can compose them without extra work. Listing 11 contains this exact statement of the composition property — we omit the proof for now, as it will be discussed in Section 6 — and illustrates how to use it to show that $f(-a) = (x \mapsto f(-x))(a)$. We will see in Section 9 how to make this even more convenient.

Listing 11: Statement and use of the composition property for the approach of Listing 2
*LeanCFC/Snippets/Intermediate.lean*

```
82 lemma cfc_comp {f g : ℂ → ℂ} {a : A} [ha : IsStarNormal a]
83     (hg : ContinuousOn g (f '' spectrum ℂ a)) (hf : ContinuousOn f (spectrum ℂ a)) :
84     cfc (g ∘ f) a = cfc g (cfc f a) := by
85   sorry
86
87 lemma cfc_comp_neg {f : ℂ → ℂ} {a : A} [ha : IsStarNormal a]
88     (hf : ContinuousOn f (- spectrum ℂ a)) :
89     cfc f (-a) = cfc (fun x ↦ f (-x)) a := by
90   rw [← Set.image_neg_eq_neg] at hf
91   change ContinuousOn f ((-id) '' spectrum ℂ a) at hf
92   rw [← cfc_id a, ← cfc_neg continuousOn_id, ← cfc_comp hf (by exact continuousOn_neg),
93       cfc_id]
94   rfl
```

## 5. Using classes to parameterize the interface

The function `cfc` from Listing 2 in the previous section addresses some of the design constraints mentioned in Section 3, especially those pertaining to usability. However, as yet we have not addressed the issues related to generality also mentioned therein. In this section, we address how to do this by providing a **`class`**-based interface to the continuous functional calculus, still in the setting of *unital* algebras. Essentially, we examine the procedure which implements Listing 2 from the isomorphism in Listing 1 and ask what features of Listing 1 are necessary to make this procedure work. As such the discussion in this section is largely orthogonal to that of Section 4.

[9] Spoilers: The left-hand side is straightforward, `a` is normal by assumption, and the result lies in the subalgebra `StarAlgebra.elemental ℂ a`. For the right-hand side, `continuousFunctionalCalculus a f` is an element of `StarAlgebra.elemental ℂ a`, and the element is normal because this algebra is commutative. Then the outer `continuousFunctionalCalculus` application (with the function `g`) is applied to this element yielding an element of the *sub*-subalgebra generated by $f(a)$ inside the subagebra generated by $a$. Then the right-hand side is coerced to the unital $C^*$-subalgebra generated by $a$, so the equality takes place in `StarAlgebra.elemental ℂ a`.

5.1. **Alternate scalar rings.** Listing 2 hard-codes the scalar ring $\mathbb{C}$, which means that utilizing functions $f : \mathbb{R} \to \mathbb{R}$ when $a$ is selfadjoint is burdensome. Even trivial considerations would become onerous if we were to use it as a general framework. For example, if $f, g : \mathbb{R} \to \mathbb{R}$, it is trivial that $(f + g)(a) = f(a) + g(a)$ because the continuous functional calculus is a $*$-homomorphism, but with the previous setup, we would first have to turn these functions into functions $f', g' : \mathbb{C} \to \mathbb{C}$. On paper, this is trivial: simply expand the codomain and recall that for selfadjoint elements the spectrum is contained in $\mathbb{R}$, so it doesn't matter how $f', g'$ are defined elsewhere.

In Lean, we would need to let `f' := fun x : ℂ ↦ (↑(f x.re : ℂ))`, and then even showing $(f + g)'(a) = f'(a) + g'(a)$ would involve multiple steps, as one first needs to establish $(f + g)' = f' + g'$ (this is not hard, but it is an extra step). The situation only gets worse when one considers composition of functions, as the function `x ↦ ↑x.re : ℂ → ℂ` is ill-behaved.

An inelegant and maintenance-heavy solution to this problem would be to simply redo everything in triplicate, once for each of the scalar semirings $\mathbb{C}$, $\mathbb{R}$ and $\mathbb{R}_{\geq 0}$. This would be hard to maintain, as whenever a given declaration is changed, its corresponding versions would also need to be updated. A much better approach is to use a class to abstract over the scalar ring. A `class` is a construct in Lean which provides a common interface to each `instance` of the class. For example, `Add` is a class in Lean providing access to the `+` notation; an instance `Add α` on the type `α` declares which binary function on `α` should be used when Lean encounters `a + b` for `a b : α`. With this approach, the mathematical content consists of providing instances of the `ContinuousFunctionalCalculus` class for $\mathbb{C}$, $\mathbb{R}$ and $\mathbb{R}_{\geq 0}$.

In our case, we will be specifying a class `ContinuousFunctionalCalculus` which provides the necessary information to construct the definition `cfc` over a given scalar semiring $R$. A naive attempt to implement this is shown in Listing 12. Instances of this class would then be provided for normal, selfadjoint or nonnegative elements of the algebra, when $R$ is $\mathbb{C}$, $\mathbb{R}$ or $\mathbb{R}_{\geq 0}$, respectively.

Listing 12: Version 1: A naive implementation of the continuous functional calculus class
*LeanCFC/Snippets/Class.lean*

```
17 class ContinuousFunctionalCalculus (R A : Type*) [CStarAlgebra A]
18     [CommSemiring R] [StarRing R] [MetricSpace R] [IsTopologicalSemiring R]
19     [ContinuousStar R] [Algebra R A] [StarModule R A] (a : A) where
20   /-- The ⋆-isomorphism underlying the continuous functional calculus for `a : A`. -/
21   toStarAlgEquiv : C(spectrum R a, R) ≃⋆ₐ[R] StarAlgebra.elemental R a
22   /-- The ⋆-isommorphism sends the identity function on `spectrum R a` to `a`. -/
23   map_id : toStarAlgEquiv (.restrict (spectrum R a) (.id R)) = ⟨a, self_mem R a⟩
```

The reader should note that we are *not* undoing all the insights learned from the previous section; our goal is instead to provide the necessary mathematical context from which we can derive `cfc` generalized to other scalar rings. Unfortunately, Listing 12 does not accomplish this, and the issues are legion. Crucially, in the case of $\mathbb{R}_{\geq 0}$, Listing 12 is not mathematically correct, as the continuous functional calculus in a unital $C^*$-algebra does not map $C(\sigma_{\mathbb{R}_{\geq 0}}(a), \mathbb{R}_{\geq 0})$ to the closed unital $*$-$\mathbb{R}_{\geq 0}$-algebra generated by $a$, but rather to the nonnegative elements in the closed unital $*$-$\mathbb{R}$-algebra generated by $a$. Indeed, let us consider the unital $C^*$-algebra $C([0, 1], \mathbb{C})$ itself with $a$ chosen to be the identity function on $[0, 1]$. In this case, $C^*_1(a)$ is the full algebra, $\sigma_{\mathbb{R}_{\geq 0}}(a) = \sigma_{\mathbb{C}}(a) = [0, 1]$, and the continuous functional calculus for $a$ is simply the identity $*$-isomorphism. Note that the closed unital $*$-$\mathbb{R}_{\geq 0}$-algebra generated by $a$ consists of uniform limits of polynomials with coefficients in $\mathbb{R}_{\geq 0}$, and therefore any function $f : [0, 1] \to \mathbb{R}_{\geq 0}$ in this subalgebra satisfies $f(x) \geq f(0)$ for all $x \in [0, 1]$. But not all functions in the algebra $C(\sigma_{\mathbb{R}_{\geq 0}}(a), \mathbb{R}_{\geq 0}) = C([0, 1], \mathbb{R}_{\geq 0})$ satisfy this property, as witnessed by the function $x \mapsto 1 - x$.

The solution to this mathematical dilemma is to no longer require that this is a $*$-*iso*morphism, and instead to require that it is a $*$-*homo*morphism into $\mathcal{A}$, but this approach calls for slightly more finesse. In order to recover some of the nice properties of the continuous functional calculus, we need to somehow "remember" properties of the $*$-isomorphism which we can no longer reference directly. In the next attempt in Listing 13, we achieve this by means of the `closedEmbedding` and `map_spectrum` conditions. The former guarantees that the function is continuous, the range is closed, and the topology on the continuous functions $C(\sigma_R(a), R)$ coincides with the pullback of the topology on $\mathcal{A}$ through the $*$-homomorphism `toStarAlgHom`. The latter is necessary because the *spectral permanence* property is not automatic over other scalar rings.

When $R$ is either $\mathbb{C}$ or $\mathbb{R}$, by the Stone–Weierstrass theorem, one can recover the $*$-isomorphism in Listing 12 from this closed embedding, which shows that we are not relinquishing any provability when $R := \mathbb{C}$.

To elaborate on this point concerning spectral permanence, note that if `f : C(spectrum ℂ a, ℂ)` and `a : A` is a normal element of a unital $C^*$-algebra, then, under the auspices of Listing 12, because `toStarAlgEquiv` is an algebra isomorphism, the spectrum of `toStarAlgEquiv f` coincides with the spectrum of `f`, which is `Set.range f`. Note that the spectrum of `toStarAlgEquiv f` is the spectrum *relative to the subalgebra* `StarAlgebra.elemental ℂ a`, which is a unital $C^*$-subalgebra of `A`. By spectral permanence, the spectrum of `toStarAlgEquiv f` is equal to the spectrum of `↑(toStarAlgEquiv f) : A`. Thus `Set.range f` coincides with `spectrum ℂ (↑(toStarAlgEquiv f) : A)`, which is the `map_spectrum` condition in the context of Listing 13 (since `toStarAlgHom` should be thought of as the composition of `toStarAlgEquiv` with the coercion from the subalgebra into `A`). However, if we do not assume that `A` is a unital $C^*$-algebra, then the spectral permanence property may fail; that is, the spectrum of an element in a subalgebra depends on the ambient algebra. Consequently, despite the fact that the `toStarAlgEquiv` is recoverable (for `R` equal to $\mathbb{C}$ or $\mathbb{R}$) from Listing 13, we would not be able to recover the `map_spectrum` condition because we would not know that `spectrum R (toStarAlgEquiv f)` is equal to `spectrum R (toStarAlgHom f)`.

Listing 13: Version 2: Avoiding the isomorphism
*LeanCFC/Snippets/Class.lean*

```
31 class ContinuousFunctionalCalculus (R A : Type*) [CStarAlgebra A]
32     [CommSemiring R] [StarRing R] [MetricSpace R] [IsTopologicalSemiring R]
33     [ContinuousStar R] [Algebra R A] (a : A) where
34   /-- The ⋆-homomorphism underlying the continuous functional calculus for ‘a : A‘. -/
35   toStarAlgHom : C(spectrum R a, R) →⋆ₐ[R] A
36   /-- The ⋆-homomorphism sends the identity function on ‘spectrum R a‘ to ‘a‘. -/
37   map_id : toStarAlgHom (.restrict (spectrum R a) (.id R)) = a
38   /-- The ⋆-homomorphism is a closed embedding. -/
39   closedEmbedding : IsClosedEmbedding toStarAlgHom
40   /-- The spectrum of the image of any function under the $*$-homomorphism is just
41   the range of that function. -/
42   map_spectrum (f : C(spectrum R a, R)) : spectrum R (toStarAlgHom f) = Set.range f
```

It turns out that there are still several issues with the class as given in Listing 13. The first is that, as was already the case in Listing 12, `a : A` is an argument to the class, which means that definitions which depend on it would be mired in rewriting. Concretely, if we implemented `cfc` with a typeclass argument `[ContinuousFunctionalCalculus R A a]`, we would suffer from the same rewriting failures as those caused in Listing 7 by the `[IsStarNormal a]` typeclass argument. Another issue is that Lean would not be able to find instances of `ContinuousFunctionalCalculus`, as they would depend on proofs that the element `a` satisfies some hypothesis (i.e., it is normal, selfadjoint or nonnegative), and since the type does not depend on these hypotheses, nor are these hypotheses themselves `class`es (except `IsStarNormal`), Lean would not be able to infer this information during type class synthesis. For more information about Lean's type class system, see [Baa22]. One way around this would be to add instances only of the form `ContinuousFunctionalCalculus R ↑a`, where `a` is a term of the appropriate subtype (that is to say, one of `{x : A // IsStarNormal x}`, `x : selfAdjoint A`, or `{x : A // 0 ≤ A}`), and `↑` is the coercion from that subtype into `A`. This approach would be burdensome as users would always have to bundle their elements into the appropriate subtype before applying the continuous functional calculus.

The solution to both of these problems is to reframe our perspective. Instead of thinking of a continuous functional calculus as associated to an element $a \in \mathcal{A}$, we should instead think of it as a property of the algebra $\mathcal{A}$ itself which only provides information for elements satisfying a certain predicate. Note that generalizing to other scalar rings comes with the additional requirement that we must also generalize the predicate on elements of the algebra that are guaranteed to have a continuous functional calculus. In particular, when the scalar ring is $\mathbb{C}$, $\mathbb{R}$ or $\mathbb{R}_{\geq 0}$ (the only ones we care about), the corresponding predicate on $a \in \mathcal{A}$ is that $a$ is normal, selfadjoint, or nonnegative, respectively. This predicate depends only on the scalar ring, and not the algebra itself. Because it doesn't change, we mark it as an `outParam` in the class definition, which means that Lean will automatically infer the value of this parameter when it is needed; this avoids the need to provide it explicitly every time we use the continuous functional calculus, which would be cumbersome.

Essentially, the way Lean infers the predicate `p` is by delaying its assignment until it finds an instance of `ContinuousFunctionalCalculus R ?q` for *some* predicate `?q : A → Prop`, and then assigning `p` to be `?q`. We therefore arrive at a workable version of the continuous functional calculus class, as shown in Listing 14.

Listing 14: Version 3: Bundling the predicate
*LeanCFC/Snippets/Class.lean*

```
49 class ContinuousFunctionalCalculus (R A : Type*) (p : outParam (A → Type*))
50     [CStarAlgebra A] [CommSemiring R] [StarRing R] [MetricSpace R]
51     [IsTopologicalSemiring R] [ContinuousStar R] [Algebra R A] where
52   toStarAlgHom {a} (ha : p a) : C(spectrum R a, R) →⋆ₐ[R] A
53   map_id {a} (ha : p a) : toStarAlgHom ha (.restrict (spectrum R a) (.id R)) = a
54   closedEmbedding {a} (ha : p a) : IsClosedEmbedding (toStarAlgHom ha)
55   map_spectrum {a} (ha : p a) (f : C(spectrum R a, R)) :
56     spectrum R (toStarAlgHom ha f) = Set.range f
57   predicate_preserving {a} (ha : p a) (f : C(spectrum R a, R)) : p (toStarAlgHom ha f)
```

Note that there is now an additional field: `predicate_preserving`. This is valuable in order to establish the composition property without extraneous hypotheses.

5.2. **Abstracting the $C^*$-algebra requirement.** With the version of the continuous functional calculus class as given in Listing 14 we still cannot create an instance of this class for `Matrix n n ℝ`, or even `Matrix n n ℂ`. This is because the class requires that the algebra $\mathcal{A}$ be a unital $C^*$-algebra. The type `Matrix n n ℂ` does not have a (globally available) $C^*$-algebra instance in Mathlib because, as discussed in the introduction, Mathlib prefers not to choose a norm on matrices. Likewise, `Matrix n n ℝ` can never have a $C^*$-algebra instance because it is not an algebra *over* $\mathbb{C}$. Fortunately, this is easily remedied: we simply remove the requirement that $\mathcal{A}$ is a unital $C^*$-algebra in the definition of the class, and replace it with much weaker requirements. Namely, we only require that $\mathcal{A}$ is a unital topological $*$-$R$-algebra. The changes to the type class assumptions are shown in Listing 15, where the structure fields coincide with those from Listing 14.

Listing 15: Version 4: Removing the $C^*$-algebra constraint
*LeanCFC/Snippets/Class.lean*

```
64 class ContinuousFunctionalCalculus (R A : Type*) (p : outParam (A → Type*))
65     [CommSemiring R] [StarRing R] [MetricSpace R] [IsTopologicalSemiring R]
66     [ContinuousStar R] [Ring A] [StarRing A] [TopologicalSpace A] [Algebra R A] where
67   toStarAlgHom {a} (ha : p a) : C(spectrum R a, R) →⋆ₐ[R] A
68   map_id {a} (ha : p a) : toStarAlgHom ha (.restrict (spectrum R a) (.id R)) = a
69   closedEmbedding {a} (ha : p a) : IsClosedEmbedding (toStarAlgHom ha)
70   map_spectrum {a} (ha : p a) (f : C(spectrum R a, R)) :
71     spectrum R (toStarAlgHom ha f) = Set.range f
72   predicate_preserving {a} (ha : p a) (f : C(spectrum R a, R)) : p (toStarAlgHom ha f)
```

5.3. **Data and propositions.** The final problem that remains with the class in Listing 15 is that it contains *data* in the form of `toStarAlgHom`. The issue is that, for type classes to work effectively, the data in instances should be unique up to *definitional*, or *judgmental*, equality. There is no problem when all the fields of a class are `Prop`-valued because in Lean any two proofs of a proposition are definitionally equal by fiat, which is called *proof irrelevance*, also known as *Uniqueness of Identity Proofs (UIP)*.

Mathematicians generally don't distinguish between definitional equality and *propositional* equality (i.e., provable equality), but in dependently typed proof assistants with an intensional type theory, like Lean, the distinction can be important. We won't delve into the details of this distinction here, but as a rough approximation: two terms are definitionally equal if, after unfolding all definitions and without applying any theorems, the terms are the same. In Lean, one can check whether `a b : α` are definitionally equal by writing `example : a = b := rfl`; `a` and `b` are definitionally equal if and only if this compiles successfully. Two terms are propositionally equal if there is a proof that they are equal. Therefore, terms which are definitionally equal are also propositionally equal, but not conversely.

In order for Listing 15 to work effectively as a class, we would need to ensure that the `toStarAlgHom` field is unique up to definitional equality for any two instances of `ContinuousFunctionalCalculus` we provide to Lean. This is actually problematic, as the following example with `A := Matrix n n ℂ` shows. As mentioned previously, Mathlib does not equip matrices with a $C^*$-algebra structure, at least not globally[10]. However, an instance of the class `ContinuousFunctionalCalculus ℂ IsStarNormal` for `Matrix n n ℂ` is still desirable, and this could be provided by means of diagonalization[11]. The problem arises if one ever activates, even temporarily, a $C^*$-algebra structure on `Matrix n n ℂ`, as then there would be two instances of the continuous functional calculus available (one via abstract $C^*$-algebra theory and the Gelfand transform, and one via diagonalization), and their underlying $*$-homomorphisms would not be definitionally equal.

Our solution relies on the following observation: the continuous functional calculus is *unique* due to the Stone–Weierstrass theorem (at least when $R$ is $\mathbb{C}$ or $\mathbb{R}$; $\mathbb{R}_{\geq 0}$ is more finicky as we discuss in Section 6). Therefore, instead of providing the $*$-homomorphisms as *data* in the class, we instead only require that there exist $*$-homomorphisms with the specified properties, thereby turning `ContinuousFunctionalCalculus` into a `Prop`-valued class and avoiding the problem of definitional equality entirely. This leads us to Listing 16, which is the version of `ContinuousFunctionalCalculus` present in Mathlib.

Listing 16: Version 5: Omitting data
*LeanCFC/Snippets/Class.lean*

```
79 class ContinuousFunctionalCalculus (R A : Type*) (p : outParam (A → Prop))
80     [CommSemiring R] [StarRing R] [MetricSpace R] [IsTopologicalSemiring R]
81     [ContinuousStar R] [Ring A] [StarRing A] [TopologicalSpace A]
82     [Algebra R A] : Prop where
83   predicate_zero : p 0
84   [compactSpace_spectrum (a : A) : CompactSpace (spectrum R a)]
85   spectrum_nonempty [Nontrivial A] (a : A) (ha : p a) : (spectrum R a).Nonempty
86   exists_cfc_of_predicate : ∀ a, p a → ∃ φ : C(spectrum R a, R) →⋆ₐ[R] A,
87     IsClosedEmbedding φ ∧ φ ((ContinuousMap.id R).restrict <| spectrum R a) = a ∧
88       (∀ f, spectrum R (φ f) = Set.range f) ∧ ∀ f, p (φ f)
```

The first three fields of this class are unrelated to our prior discussion. In fact, these fields are not strictly necessary, but they allow us to avoid adding type class assumptions downstream when developing general theory pertaining to the continuous functional calculus and are therefore convenient. In our first iteration that was merged to Mathlib, they were not present, and we will not mention them further. It is the last field `exists_cfc_of_predicate` which contains the key statement, now bundled into a singled proposition.

Given an instance of `ContinuousFunctionalCalculus R p`, we extract a $*$-homomorphism `φ` appearing in Listing 16 into its own definiton:

Listing 17: Extracting the $*$-homomorphism from the class
*LeanCFC/Snippets/Class.lean*

```
94 noncomputable def cfcHom {a : A} (ha : p a) : C(spectrum R a, R) →⋆ₐ[R] A :=
95   (ContinuousFunctionalCalculus.exists_cfc_of_predicate a ha).choose
```

We can also define `cfc` analogously to Listing 2.

Listing 18: Defining `cfc` using the class
*LeanCFC/Snippets/Class.lean*

```
97 open scoped Classical in
98 noncomputable irreducible_def cfc (f : R → R) (a : A) : A :=
```

[10] What we mean here is that such instances are *scoped*, meaning they are only activated within a specified namespace. For example, in order to activate a $C^*$-algebra structure on matrices in Mathlib, one must `open Matrix.L2OpNorm`.

[11] To be a bit more specific, if `a : Matrix n n ℂ` is normal, then there is a unitary matrix `u : Matrix n n ℂ` and a diagonal matrix `d : Matrix n n ℂ` such that `a = uᴴ * d * u`, where `uᴴ` denotes the hermitian conjugate. Then the map `fun f : ℂ → ℂ ↦ uᴴ * (d.map f) * u` is the desired $*$-homomorphism underlying the continuous functional calculus.

```
99   if h : p a ∧ ContinuousOn f (spectrum R a)
100    then cfcHom h.1 ⟨_, h.2.restrict⟩
101    else 0
```

## 6. Uniqueness

While the fact that the continuous functional calculus is a $*$-homomorphism is important, its most salient feature is the composition property. With the final version of the continuous functional calculus in hand, we can now state this property, the same as that appearing in Listing 11 but now letting `R` and `p` vary, concisely in Lean:

Listing 19: The composition property in Lean
*LeanCFC/Snippets/Uniqueness.lean*

```
6  lemma cfc_comp {R A : Type*} (p : outParam (A → Prop))
7      [CommSemiring R] [StarRing R] [MetricSpace R] [IsTopologicalSemiring R]
8      [ContinuousStar R] [Ring A] [StarRing A] [TopologicalSpace A]
9      [Algebra R A] [ContinuousFunctionalCalculus R A p]
10     [ContinuousMap.UniqueHom R A]
11     (g : R → R) (f : R → R) (a : A) (ha : p a)
12     (hg : ContinuousOn g (f '' spectrum R a))
13     (hf : ContinuousOn f (spectrum R a)) :
14     cfc (g ∘ f) a = cfc g (cfc f a) := by
```

In the literature, this property is sometimes not stated *at all* ([Lin01, Corollary 1.3.6] or [KR97, Theorem 4.1.3 for selfadjoint elements], stated *without proof* [Tak01, p. 19], proven for polynomials and then a direct appeal to the Stone–Weierstrass theorem [Dav96, Corollary I.3.3], or by appeal to the *uniqueness* of the continuous functional calculus [KR97, Theorems 4.4.5 and 4.4.8 for normal elements] or [Bou19, Proposition I.6.7 and Corollary I.6.2]).

We sketch the latter two approaches below, after recalling some preliminaries. For the sake of clarity, we will temporarily denote by $\mathrm{cfc}(f, x)$ the application of the continuous functional calculus for $x$, a normal element in a unital $C^*$-algebra, to a function $f \in C(\sigma_{\mathbb{C}}(x), \mathbb{C})$. Recall that this is a $*$-homomorphism of the first argument and additionally satisfies $\mathrm{cfc}(\mathrm{id}, a) = a$.

Given any complex polynomial $p$ in two variables, and $x$ a normal element in a unital $C^*$-algebra $A$, it is standard practice to evaluate $p$ at $x$ by $p_A(x) := p(x, x^*)$. It is straightforward to see that evaluation of $p$ at a normal element commutes with $*$-homomorphisms. That is, if $\phi : A \to B$ is a $*$-homomorphism and $x \in A$ is normal, then $\phi(p_A(x)) = p_B(\phi(x))$.

On the other hand, $p_{\mathbb{C}} : z \mapsto p(z, \bar{z})$ is a continuous complex-valued function defined on $\mathbb{C}$, which means that an alternative way to evaluate $p$ at the normal element $x$ is given by $\mathrm{cfc}(p_{\mathbb{C}}, x)$. It is straightforward to check that these two notions of evaluation coincide :

$$p_A(x) = p_A(\mathrm{cfc}(\mathrm{id}, x)) = \mathrm{cfc}(p_{C(\sigma_{\mathbb{C}}(x), \mathbb{C})}(\mathrm{id}), x) = \mathrm{cfc}(p_{\mathbb{C}}, x).$$

*Proof via polynomials.* Let $a \in \mathcal{A}$ be a normal element in a unital $C^*$-algebra, and let $f : \sigma_{\mathbb{C}}(a) \to \mathbb{C}$ be a continuous function. Given an arbitrary complex polynomial $p$ in two variables, we use the commutation of polynomials with $*$-homomorphisms to obtain

$$\begin{aligned}\mathrm{cfc}(p_{\mathbb{C}} \circ f, a) &= \mathrm{cfc}(p_{C(\sigma_{\mathbb{C}}(a), \mathbb{C})}(f), a) \\ &= p_A(\mathrm{cfc}(f, a)) \\ &= \mathrm{cfc}(p_{\mathbb{C}}, \mathrm{cfc}(f, a)).\end{aligned}$$

By the Stone–Weierstrass theorem, the set of all such polynomial functions $p_{\mathbb{C}}$ is dense in the unital $C^*$-algebra $C(\sigma_{\mathbb{C}}(\mathrm{cfc}(f, a)), \mathbb{C})$. Moreover, since $*$-homomorphisms between $C^*$-algebras are continuous, we obtain by taking limits $\mathrm{cfc}(g \circ f, a) = \mathrm{cfc}(g, \mathrm{cfc}(f, a))$ for any continuous $g : \sigma_{\mathbb{C}}(\mathrm{cfc}(f, a)) \to \mathbb{C}$. □

*Proof via uniqueness.* Let $a \in \mathcal{A}$ be a normal element in a unital $C^*$-algebra. Note that if $\phi : C(\sigma_{\mathbb{C}}(a), \mathbb{C}) \to \mathcal{A}$ is any $*$-homomorphism such that $\phi(\mathrm{id}) = a$, then $\phi(f) = \mathrm{cfc}(f, a)$ for any continuous $f : \sigma_{\mathbb{C}}(a) \to \mathbb{C}$.

Indeed, for any polynomial $p$ in two variables we have

$$\phi(p_{\mathbb{C}}) = \phi(p_{C(\sigma_{\mathbb{C}}(a),\mathbb{C})}(\mathrm{id})) = p_A(\phi(\mathrm{id})) = p_A(a) = \mathrm{cfc}(p_{\mathbb{C}}, a),$$

then by the Stone–Weierstrass theorem and the continuity of $\phi$, this equality extends to all continuous functions. Therefore the continuous functional calculus is the *unique* $*$-homomorphism from $C(\sigma_{\mathbb{C}}(a), \mathbb{C})$ to $\mathcal{A}$ sending the identity function to $a$.

Now suppose that $f : \sigma_{\mathbb{C}}(a) \to \mathbb{C}$ is continuous. There is the continuous functional calculus for the normal element $\mathrm{cfc}(f, a)$, which sends a continuous function $g : \sigma_{\mathbb{C}}(\mathrm{cfc}(f, a)) \to \mathbb{C}$ to $\mathrm{cfc}(g, f(a))$. However, for such $g$, notice that $g \circ f$ is a continuous function on $\sigma_{\mathbb{C}}(a)$, and so $\mathrm{cfc}(g \circ f, a)$ is also meaningful, but is *a priori* different from $\mathrm{cfc}(g, \mathrm{cfc}(f, a))$. But the mapping $g \mapsto \mathrm{cfc}(g \circ f, a)$ is a $*$-homomorphism from $C(\sigma_{\mathbb{C}}(\mathrm{cfc}(f, a)), \mathbb{C})$ to $\mathcal{A}$ (it's the continuous functional calculus for $a$ composed with the $*$-homomorphism given by precomposition with $f$). Finally, $\mathrm{cfc}(\mathrm{id} \circ f, a) = \mathrm{cfc}(f, a)$, and so by uniqueness, we must have $\mathrm{cfc}(g \circ f, a) = \mathrm{cfc}(g, \mathrm{cfc}(f, a))$ for any continuous function $g : \sigma_{\mathbb{C}}(\mathrm{cfc}(f, a)) \to \mathbb{C}$. ☐

While both proofs appear in the literature, arguments akin to the former are ubiquitous throughout. For example, even though Kadison–Ringrose [KR97, Theorem 4.4.8] uses the uniqueness proof for the composition property, in [KR97, Propostion 4.2.3] the authors rely on the polynomial argument to show that the positive and negative parts of a selfadjoint element are unique, where a uniqueness argument would have worked just as well. Indeed, [KR97, Propostion 4.2.3] is implemented in Mathlib using uniqueness (`CFC.posPart_negPart_unique`).

Although the polynomial argument occurs frequently in the literature, we instead adopt the uniqueness approach in Mathlib for two primary reasons. The first is practical: limiting arguments by appeal to hand-waving are acceptable on paper, but formalizing them can be cumbersome; uniqueness is simply easier. Secondly, we also need Listing 19 to hold when $R := \mathbb{R}_{\geq 0}$ in order for the continuous functional calculus to be useful, and unfortunately the Stone–Weierstrass theorem does not hold in this case. Moreover, as we indicated previously, not every nonnegative continuous function with domain $s \subseteq \mathbb{R}_{\geq 0}$ is a uniform limit of polynomials *with nonnegative coefficients* since such limits must be monotone functions on the nonnegative real numbers because each monomial is monotone. The reader may have noticed that the uniqueness argument *also* makes use of the Stone–Weierstrass theorem. Despite this, it is still possible to prove uniqueness of continuous unital $\mathbb{R}_{\geq 0}$-$*$-algebra homomorphisms with domain $C(s, \mathbb{R}_{\geq 0})$ for $s$ compact (and hence the composition property for the continuous functional calculus over $\mathbb{R}_{\geq 0}$), provided the codomain is a unital topological $\mathbb{R}$-$*$-algebra. So, we define the following **`class`** in Mathlib:

Listing 20: Continuous Functional Calculus Uniqueness Class
*LeanCFC/Snippets/Uniqueness.lean*

```
28 class ContinuousMap.UniqueHom (R A : Type*) [CommSemiring R] [StarRing R]
29     [MetricSpace R] [IsTopologicalSemiring R] [ContinuousStar R] [Ring A] [StarRing A]
30     [TopologicalSpace A] [Algebra R A] : Prop where
31   eq_of_continuous_of_map_id (s : Set R) [CompactSpace s]
32     (φ ψ : C(s, R) →⋆ₐ[R] A) (hφ : Continuous φ) (hψ : Continuous ψ)
33     (h : φ (.restrict s <| .id R) = ψ (.restrict s <| .id R)) :
34     φ = ψ
```

Having this as a class allows us to assume generic instances of it when developing general theory for the continuous functional calculus, especially for lemmas involving the composition property. For example, the lemmas contained in the `Comp` section of *LeanCFC/Implementation/Unital.lean* all need this hypothesis, including `cfc_comp_neg` and `cfc_comp_star` which constitute Item (3) from Section 3.1.

Note that whenever $\mathcal{A}$ is a unital topological $\mathbb{K}$-algebra (with $\mathbb{K} := \mathbb{R}$ or $\mathbb{K} := \mathbb{C}$), then the Stone–Weierstrass theorem allows us to create an instance of this class. However, we can also provide an instance of this class for $R := \mathbb{R}_{\geq 0}$ provided that $\mathcal{A}$ is a topological $\mathbb{R}$-algebra by means of the following argument. The idea is, for $s \subseteq \mathbb{R}_{\geq 0}$ compact, to take a unital $*$-algebra homomorphism $\phi : C(s, \mathbb{R}_{\geq 0}) \to \mathcal{A}$ (over $\mathbb{R}_{\geq 0}$), and create a unital $*$-algebra homomorphism $\hat{\phi} : C(s, \mathbb{R}) \to \mathcal{A}$ (over $\mathbb{R}$) by defining $\hat{\phi}(f) := \phi(f_+) - \phi(f_-)$. Then the uniqueness for $\mathbb{R}$ implies uniqueness for $\mathbb{R}_{\geq 0}$. So, making uniqueness of the continuous functional calculus into a *class* allows us to unify the framework for $\mathbb{C}$, $\mathbb{R}$ and $\mathbb{R}_{\geq 0}$.

One important question to consider is: why is `ContinuousMap.UniqueHom` a separate class instead of being bundled into the `ContinuousFunctionalCalculus` class itself? For starters, this class simply has instances in much greater generality (i.e., there is no need for $\mathcal{A}$ to be a $C^*$-algebra at all). It is also dictated in part by the way we chose to obtain instances of continuous functional calculi for $\mathbb{R}_{\geq 0}$ and $\mathbb{R}$ from those of $\mathbb{R}$ and $\mathbb{C}$, respectively. We create a generic lemma for `ContinuousFunctionalCalculus` that allows us to pass from a calculus over a scalar ring $R$ with predicate $p$ to a calculus over a scalar subring $S$ with predicate $q$, provided that any $a \in \mathcal{A}$ satisfies $q$ if and only if it satisfies $p$ and the $R$-spectrum of $a$ is contained in $S$. As an example, $a \in \mathcal{A}$ is selfadjoint if and only if it is normal and its spectrum is contained in $\mathbb{R}$. Therefore we obtain a continuous functional calculus over $\mathbb{R}$ for selfadjoint elements from the one over $\mathbb{C}$ for normal elements. We refer to Section 8 for a more precise description of this restriction process. If we were to include the uniqueness criterion in the main class, then we would be forced to prove it when we establish this generic restriction lemma; this would be burdensome or impossible in that generality.

Another more important reason: we don't want uniqueness to hold only for $*$-homomorphisms into the algebra on which we have the continuous functional calculus for a given element, but for any suitable algebra. This is essential to prove, for example, that the continuous functional calculus commutes with $*$-homomorphisms between unital $C^*$-algebras. To sketch the idea, let $\mathcal{A}$ and $\mathcal{B}$ be unital $C^*$-algebras and $\phi : \mathcal{A} \to \mathcal{B}$ a $*$-homomorphism and $a \in \mathcal{A}$ a normal element. To show that $\phi(f(a)) = f(\phi(a))$ for any continuous function $f : \sigma_{\mathbb{C}}(a) \to \mathbb{C}$ we apply uniqueness of $*$-homormorphisms from $C(\sigma_{\mathbb{C}}(a), \mathbb{C})$ to $\mathcal{B}$ (note that $a \in \mathcal{A}$ but the codomain is $\mathcal{B}$!) taking a specified value on the (restriction of) the identity function. Note $f \mapsto \phi(f(a))$ is a composition of $*$-homomorphisms ($\phi$ and the continuous functional calculus for $a$) which sends the identity function to $\phi(a)$. Likewise, $f \mapsto f(\phi(a))$ is also a composition of $*$-homomorphisms (the continuous functional calculus for $\phi(a)$ and the restriction map from $C(\sigma_{\mathbb{C}}(a), \mathbb{C})$ to $C(\sigma_{\mathbb{C}}(\phi(a)), \mathbb{C})$) sending the identity function to $\phi(a)$. By uniqueness, these must coincide.

## 7. The continuous functional calculus for non-unital algebras

Non-unital $C^*$-algebras play a prominent role in the subject but, of course, the spectrum does not make sense *a priori* for non-unital algebras, as there is no notion of invertibility. In textbooks on the subject, the approach is as follows: given a non-unital $C^*$-algebra $\mathcal{A}$, construct the *minimal unitization* $\mathcal{A}^{+1}$. This is the unital $C^*$-algebra $\mathcal{A}^{+1} := \mathbb{C} \times \mathcal{A}$ where multiplication is defined by $(z, a)(w, b) := (zw, zb + wa + ab)$ and the norm is defined by $\|(z, a)\| := \max\{|z|, \|x \mapsto zx + ax\|\}$. It is called the minimal unitization because, for any unital $C^*$-algebra $\mathcal{B}$ and any $*$-homomorphism $\phi : \mathcal{A} \to \mathcal{B}$, there is a unique *unital* $*$-homomorphism $\hat{\phi} : \mathcal{A}^{+1} \to \mathcal{B}$ extending $\phi$[12].

At this point, textbook approaches diverge on their treatment, although all of these approaches are ultimately equivalent. We list a few of these approaches here. Kadison and Ringrose [KR97], for example, never touch non-unital $C^*$-algebras, as their focus primarily drifts to von Neumann algebras, which are always unital.

Lin [Lin01], Pedersen [Ped79] and Blackadar [Bla06] take the approach of defining the spectrum of $a \in \mathcal{A}$ as the spectrum of $(0, a) \in \mathcal{A}^{+1}$, but only when $\mathcal{A}$ does not contain any identity element (recall that, for us, a non-unital algebra may still contain an identity). When $\mathcal{A}$ is unital, the spectrum has the usual meaning[13]. This approach is also employed by Davdison [Dav96] and Fillmore [Fil96], but they never explicitly define the spectrum for non-unital algebras, and instead freely pass to the unitization whenever it is convenient. In these approaches, the continuous functional calculus is a $*$-homomorphism from $C_0(\sigma_{\mathbb{C}}(a) \setminus \{0\}, \mathbb{C})$ to $\mathcal{A}$. Here, $C_0$ denotes the non-unital algebra of continuous functions *vanishing at infinity*, i.e., $f \in C_0(\sigma_{\mathbb{C}}(a) \setminus \{0\}, \mathbb{C})$ if, for every $\epsilon > 0$, $\{z \in \sigma_{\mathbb{C}}(a) \setminus \{0\} \mid |f(z)| \geq \epsilon\}$ is compact. The reader should take note: this does not mean $f(z) \to 0$ as $|z| \to \infty$, as this wouldn't even make sense for $f \in C_0(\sigma_{\mathbb{C}}(a) \setminus \{0\}, \mathbb{C})$. Rather, since $\sigma_{\mathbb{C}}(a)$ is compact, then $f \in C(\sigma_{\mathbb{C}}(a) \setminus \{0\}, \mathbb{C})$ vanishes at infinity when $f(z) \to 0$ as $z \to 0$, with this condition being vacuous when $0 \notin \sigma_{\mathbb{C}}(a)$.

[12]There is another widely used notion of "minimal unitization", namely the algebra $\widetilde{A}$ which is just $\mathcal{A}$ if $\mathcal{A}$ happens to be unital, and $\mathcal{A}^{+1}$ if $\mathcal{A}$ is not unital. By contrast, $\mathcal{A}^{+1}$ always adds a *new* identity element. This variant $\widetilde{\mathcal{A}}$, which we shall not use in the rest of the paper, satisfies a different universal property: for any unital $C^*$-algebra $\mathcal{B}$ and any *nondegenerate* $*$-homomorphism $\phi : \mathcal{A} \to \mathcal{B}$, there is a unique *unital* $*$-homomorphism $\hat{\phi} : \widetilde{\mathcal{A}} \to \mathcal{B}$ extending $\phi$. Here, a $*$-homomomorphism is said to be *nondegenerate* if it maps any approximate unit to an approximate unit.

[13]In other words, they define the spectrum of $a \in \mathcal{A}$ as the spectrum of $a$ in the unitization $\widetilde{\mathcal{A}}$ described in the previous footnote.

Bourbaki [Bou19] and Dixmier [Dix69] take a different approach by defining two notions of spectrum: the usual spectrum in unital algebras, and another one for non-unital algebras which is the spectrum in the unitization. The latter of these always contains zero. Bourbaki and Dixmier provide alternative notation for these two notions of spectrum, but Bourbaki at least still refers to the latter as "spectre de $x$ relativement à $\mathcal{A}$". In this variation, the non-unital continuous functional calculus is a $*$-homomorphism from the ideal of functions $f \in C(\sigma_{\mathbb{C}}((0,a)), \mathbb{C})$ for which $f(0) = 0$ to $\mathcal{A}$. We denote this ideal as $C(\sigma_{\mathbb{C}}((0,a)), \mathbb{C})_0$. Please be careful not to confuse this with the continuous functions vanishing at infinity, noting the placement of the subscript 0. Takesaki [Tak01, p. 7] utilizes the same technique, but refers to the spectrum in the unitization as the *quasispectrum*.

From the perspective of formalization, all of these textbook approaches have the downside that they depend upon the construction of another type — the unitization (i.e., they are *extrinsic* formulations). In the case of [Lin01, Ped79, Bla06], the definition is even more problematic as it requires a case split on whether the algebra is unital or not. Our preferred approach is closest to Takesaki's in [Tak01], but we opt for an *intrinsic* definition of the quasispectrum for non-unital algebras.

Given a non-unital algebra $\mathcal{A}$ over a semifield $R$, an element $r \in R$ is in the *quasispectrum* of $a \in \mathcal{A}$ if either $r = 0$ or $x := -(r^{-1}a)$ is not *quasiregular* (i.e., there is no $y \in \mathcal{A}$ such that $x + y + xy = 0 = y + x + yx$). We note that $x$ is quasiregular in $\mathcal{A}$ if and only if $(1, x) = r^{-1} \cdot (r, -a)$ is invertible in $\mathcal{A}^{+1}$; consequently the quasispectrum of $a$ in $\mathcal{A}$ coincides with the spectrum of $(0, a)$ in $\mathcal{A}^{+1}$. It therefore coincides with the quasispectrum of $a$ in the sense of [Tak01, p. 7]. We generally use the notation $\sigma_{n,R}(a)$ (or `σₙ R a` in Lean) for the quasispectrum of $a \in \mathcal{A}$ in the non-unital $R$-algebra $\mathcal{A}$.

Listing 21: Definition of the quasispectrum
*LeanCFC/Snippets/NonUnital.lean*

```
5  /-- If `A` is a non-unital `R`-algebra, the `R`-quasispectrum of `a : A` consists of those
6  `r : R` such that if `r` is invertible (in `R`), then `-(r⁻¹ • a)` is not quasiregular.
7
8  The quasispectrum is precisely the spectrum in the unitization when `R` is a
9  commutative ring. -/
10 def quasispectrum (R : Type*) {A : Type*} [CommSemiring R] [NonUnitalRing A] [Module R A]
11     (a : A) : Set R :=
12   {r : R | (hr : IsUnit r) → ¬ IsQuasiregular (-(hr.unit⁻¹ • a))}
```

This leads us to a definition of the continuous functional calculus for non-unital algebras which closely resembles the unital one. In Listing 22, `C(σₙ R a, R)₀` denotes the type of continuous functions which map zero to zero.

Listing 22: Non-unital continuous functional calculus
*LeanCFC/Snippets/NonUnital.lean*

```
17 class NonUnitalContinuousFunctionalCalculus (R : Type*) {A : Type*} (p : outParam (A → Prop))
18     [CommSemiring R] [Nontrivial R] [StarRing R] [MetricSpace R] [IsTopologicalSemiring R]
19     [ContinuousStar R] [NonUnitalRing A] [StarRing A] [TopologicalSpace A] [Module R A]
20     [IsScalarTower R A A] [SMulCommClass R A A] : Prop where
21   predicate_zero : p 0
22   [compactSpace_quasispectrum : ∀ a : A, CompactSpace (σₙ R a)]
23   exists_cfc_of_predicate : ∀ a, p a → ∃ φ : C(σₙ R a, R)₀ →⋆ₙₐ[R] A,
24     IsClosedEmbedding φ ∧ φ ⟨(ContinuousMap.id R).restrict <| σₙ R a, rfl⟩ = a ∧
25       (∀ f, σₙ R (φ f) = Set.range f) ∧ ∀ f, p (φ f)
```

As in the unital case, we also provide a uniqueness class for the non-unital continuous functional calculus.

Listing 23: Non-Unital Continuous Functional Calculus Uniqueness Class
*LeanCFC/Snippets/NonUnital.lean*

```
27 class ContinuousMapZero.UniqueHom (R A : Type*) [CommSemiring R] [StarRing R]
28     [MetricSpace R] [IsTopologicalSemiring R] [ContinuousStar R] [NonUnitalRing A] [StarRing A]
29     [TopologicalSpace A] [Module R A] [IsScalarTower R A A] [SMulCommClass R A A] : Prop where
30   eq_of_continuous_of_map_id (s : Set R) [CompactSpace s] [Zero s] (h0 : (0 : s) = (0 : R))
31     (φ ψ : C(s, R)₀ →⋆ₙₐ[R] A) (hφ : Continuous φ) (hψ : Continuous ψ)
32     (h : φ (⟨.restrict s <| .id R, h0⟩) = ψ (⟨.restrict s <| .id R, h0⟩)) :
33     φ = ψ
```

Likewise we define `cfcₙHom` and `cfcₙ`.

Listing 24: Extracting the non-unital $*$-homomorphism from the class
*LeanCFC/Snippets/NonUnital.lean*

```
35 variable {R A : Type*} {p : A → Prop} [CommSemiring R] [Nontrivial R] [StarRing R]
36     [MetricSpace R] [IsTopologicalSemiring R] [ContinuousStar R] [NonUnitalRing A] [StarRing A]
37     [TopologicalSpace A] [Module R A] [IsScalarTower R A A] [SMulCommClass R A A]
38     [NonUnitalContinuousFunctionalCalculus R p]
39
40 noncomputable def cfcₙHom {a : A} (ha : p a) : C((σₙ R a), R)₀ →⋆ₙₐ[R] A :=
41   (NonUnitalContinuousFunctionalCalculus.exists_cfc_of_predicate a ha).choose
```

Listing 25: Defining `cfcₙ` using the class
*LeanCFC/Snippets/NonUnital.lean*

```
43 open scoped Classical in
44 noncomputable irreducible_def cfcₙ (f : R → R) (a : A) : A :=
45   if h : p a ∧ ContinuousOn f (σₙ R a) ∧ f 0 = 0
46     then cfcₙHom h.1 ⟨⟨_, h.2.1.restrict⟩, h.2.2⟩
47     else 0
```

## 8. Instantiating the continuous functional calculus

From the beginning of Section 5 up to this point, we have been discussing the *design* of the continuous functional calculus, but not how to actually *instantiate* it. A **`class`** in Lean is in some sense an abstract interface, but without **`instance`**s of the class, it serves no purpose. The mathematical content comes therefore in providing instances of the continuous functional calculus. Of course, Listing 1 leads straightforwardly to an instance of `ContinuousFunctionalCalculus ℂ IsStarNormal` for unital $C^*$-algebras. But there are several other instances we need to provide, which include versions over $\mathbb{R}$ and $\mathbb{R}_{\geq 0}$ for unital algebras, as well as non-unital versions of all these. These various instances and their interactions are summarized in Figure 1.

8.1. **Passing to scalar subrings.** On paper, other than proving that selfadjoint elements have real spectrum, there is virtually nothing to be be done to obtain the continuous functional calculus over $\mathbb{R}$ for selfadjoint elements from the one over $\mathbb{C}$ for normal elements. This is so much the case that it is often not even mentioned. For nonnegative elements, generally even less is said, as in textbooks nonnegative elements are often defined as selfadjoint elements with nonnegative spectrum.

In contrast, we definitely need this in Lean because $\mathbb{R}$ and $\mathbb{C}$ are separate types, and it is much nicer, where possible, to work with $\mathbb{R}$ directly rather than with its image in $\mathbb{C}$. But at the same time, this needs to be explicit. Because we have to do this four times ($\mathbb{C}$ to $\mathbb{R}$, $\mathbb{R}$ to $\mathbb{R}_{\geq 0}$, for both unital and non-unital functional calculi), we developed a suitable interface for indicating that the spectrum of an element with respect to a given scalar semiring can be reinterpreted as the spectrum relative to a scalar subring. We introduce the structure `QuasispectrumRestricts`, described in Listing 26, for this purpose. In fact, we use the exact same structure (up to definitional equality) for the *spectrum* in unital algebras as well. This is because, over a (semi)field, the quasispectrum is contained in a scalar subfield if and only if the spectrum is contained in that subfield, since they differ at most by the presence of zero.

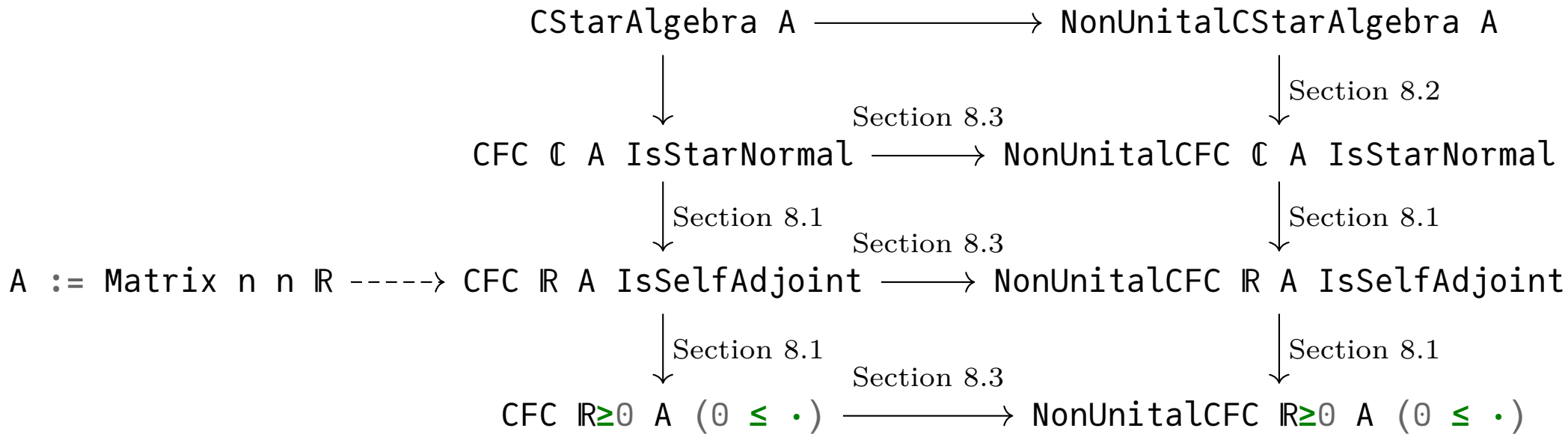

FIGURE 1. Instance graph for the continuous functional calculus. We abbreviate `ContinuousFunctionalCalculus` to `CFC` for readability. Each full arrow `X` → `Y` in this graph corresponds to an instance of the form `instance ... [X] : Y`, with annotations referencing the relevant subsection. Finally, the dashed arrow represents a specific instance of type `CFC ℝ (Matrix n n ℝ) IsSelfAdjoint`.

Listing 26: Restriction of the (quasi)spectrum
*LeanCFC/Snippets/Instances.lean*

```
 5 /-- Given an element `a : A` of an `S`-algebra, where `S` is itself an `R`-algebra, we say
 6 that the spectrum of `a` restricts via a function `f : S → R` if `f` is a left inverse of
 7 `algebraMap R S`, and `f` is a right inverse of `algebraMap R S` on `spectrum S a`.
 8
 9 For example, when `f = Complex.re` (so `S := ℂ` and `R := ℝ`), `SpectrumRestricts a f`
10 means that the `ℂ`-spectrum of `a` is contained within `ℝ`. This arises naturally when
11 `a` is selfadjoint and `A` is a C⋆-algebra.  -/
12 structure QuasispectrumRestricts {R S A : Type*} [CommSemiring R] [CommSemiring S]
13     [NonUnitalRing A] [Module R A] [Module S A] [Algebra R S] (a : A) (f : S → R) : Prop where
14   /-- `f` is a right inverse of `algebraMap R S` when restricted to `quasispectrum S a`. -/
15   rightInvOn : (quasispectrum S a).RightInvOn f (algebraMap R S)
16   /-- `f` is a left inverse of `algebraMap R S`. -/
17   left_inv : Function.LeftInverse f (algebraMap R S)
```

We then prove a general theorem (Listing 27) that if we have a continuous functional calculus over a scalar ring $S$ for a unital algebra $\mathcal{A}$ with predicate $q$, and $p$ is some other predicate which is equivalent to "$q$ and the spectrum restricts to a scalar subring $R$", then we can obtain a continuous functional calculus over $R$ for $\mathcal{A}$ with predicate $p$, under the assumption that the natural map from $R$ to $S$ is a uniform embedding. Of course, we can do the same for the non-unital continuous functional calculus as well over non-unital algebras. This general theorem can't be given as an *instance* because Lean would not be able to infer it in practice, but it allows us to provide very short proofs of the relevant instances (see for example, the proofs of `IsSelfAdjoint.instNonUnitalContinuousFunctionalCalculus`, `IsSelfAdjoint.instContinuousFunctionalCalculus`, `Nonneg.instNonUnitalContinuousFunctionalCalculus` or `Nonneg.instContinuousFunctionalCalculus`).

Listing 27: Restricting the continuous functional calculus to scalar subrings
*LeanCFC/Snippets/Instances.lean*

```
19 /-- Given a `ContinuousFunctionalCalculus S q`, if we form the predicate `p` for `a : A`
20 characterized by: `q a` and the spectrum of `a` restricts to the scalar subring `R` via
21 `f : C(S, R)`, then we can get a restricted functional calculus
22 `ContinuousFunctionalCalculus R p`. -/
23 theorem SpectrumRestricts.cfc {R S A : Type*} {p q : A → Prop} [Semifield R] [StarRing R]
24     [MetricSpace R] [IsTopologicalSemiring R] [ContinuousStar R] [Semifield S] [StarRing S]
25     [MetricSpace S] [IsTopologicalSemiring S] [ContinuousStar S] [Ring A] [StarRing A]
26     [Algebra S A] [Algebra R S] [Algebra R A] [IsScalarTower R S A] [StarModule R S]
27     [ContinuousSMul R S] [TopologicalSpace A] [ContinuousFunctionalCalculus S A q]
28     [CompleteSpace R] (f : C(S, R)) (halg : IsUniformEmbedding ⇑(algebraMap R S))
```

```
29    (h0 : p 0) (h : ∀ (a : A), p a ↔ q a ∧ SpectrumRestricts a ⇑f) :
30    ContinuousFunctionalCalculus R A p :=
```

8.2. **The non-unital instance for non-unital $C^*$-algebras.** We now describe how to construct an instance of the non-unital functional calculus (over $\mathbb{C}$, other scalar subrings are obtained from this via the generic framework above) for non-unital $C^*$-algebras. The most standard approach is to rely on *non-unital Gelfand duality*, which is a contravariant equivalence of categories between non-unital commutative $C^*$-algebras and pointed compact Hausdorff spaces. Note that this duality is constructed by reducing to the unital case, so every approach relies on the unital setting in some way. Given that we are not as interested in the natural isomorphisms involved in this equivalence of categories, but instead we primarily want an instance of `NonUnitalContinuousFunctionalCalculus`, we opt for a streamlined approach which differs slightly from that in most textbooks on the subject.

Suppose that $\mathcal{A}$ is a non-unital $C^*$-algebra. Then we can form the minimal unitization (over $\mathbb{C}$) $\mathcal{A}^{+1}$. Because this is a unital $C^*$-algebra, we have a continuous functional calculus on $\mathcal{A}^{+1}$ over $\mathbb{C}$ for normal elements. Then we can construct a non-unital $*$-homomorphism $\phi : C(\sigma_{n,\mathbb{C}}(a), \mathbb{C})_0 \to \mathcal{A}^{+1}$ as the composition

$$C(\sigma_{n,\mathbb{C}}(a), \mathbb{C})_0 \xrightarrow{(\uparrow)} C(\sigma_{n,\mathbb{C}}(a), \mathbb{C}) \xrightarrow{\sim} C(\sigma_{\mathbb{C}}((0,a)), \mathbb{C}) \longrightarrow \mathcal{A}^{+1}$$

where the first map (implemented as a coercion) sends any function to itself meanwhile forgetting the fact that it maps zero to zero, the second map is the $*$-isomorphism from arising from the set-equality $\sigma_{n,\mathbb{C}}(a) = \sigma_{\mathbb{C}}((0,a))$ (where $(0,a) \in \mathcal{A}^{+1}$), and the last map is the (unital!) continuous functional calculus in the unital $C^*$-algebra $\mathcal{A}^{+1}$.

This is not quite the map we want because it takes values in $\mathcal{A}^{+1}$, not in $\mathcal{A}$. The last step is then to realize that, because $C(\sigma_{n,\mathbb{C}}(a), \mathbb{C})_0$ is generated (as a non-unital topological star algebra) by the identity function (this is Stone–Weierstrass), and the image of the identity is $(0,a) \in \mathcal{A}^{+1}$, then the range of this non-unital $*$-homomorphism is actually contained in the image of $\mathcal{A}$ in $\mathcal{A}^{+1}$ (see `cfcₙAux_mem_range_inr`). Since the inclusion $\mathcal{A} \hookrightarrow \mathcal{A}^{+1}$ is a closed embedding, we can pull back the map described above to get a closed embedding (and a non-unital $*$-homomorphism) $C(\sigma_{n,\mathbb{C}}(a), \mathbb{C})_0 \to \mathcal{A}$ (`RCLike.nonUnitalContinuousFunctionalCalculus`).

8.3. **The non-unital instance for unital algebras.** At this point, we have described all of the vertical arrows in the instance graph of Figure 1. Hence, we can use the three variants of the non-unital calculus in any non-unital $C^*$-algebra $\mathcal{A}$, and if $\mathcal{A}$ is unital we can also refer to the three unital variants. Yet, we also want instances for the horizontal arrows. The first motivation is convenience: when setting-up the generic API for the functional calculus, lemmas involving both the unital and non-unital variants would otherwise need to assume the existence (and uniqueness) of both variations. More crucially, it is clear in Figure 1 that the horizontal arrows are necessary in order to obtain an instance of `NonUnitalContinuousFunctionalCalculus ℝ≥0 (Matrix n n ℝ) (0 ≤ ·)`, which appears, for example, when applying the general construction of square roots to the case of matrices.

Implementing these instances might seem so obvious as to be trivial, but there is some minor subtlety. Of course, given a $*$-homomorphism $\phi : C(\sigma_R(a), R) \to \mathcal{A}$, we want to restrict it to a non-unital $*$-homomorphism $\phi_0 : C(\sigma_{n,R}(a), R)_0 \to \mathcal{A}$. Formally, $\phi_0$ is obtained as the composition

$$C(\sigma_{n,R}(a), R)_0 \xrightarrow{(\uparrow)} C(\sigma_{n,R}(a), R) \xrightarrow{\psi} C(\sigma_R(a), R) \xrightarrow{\phi} \mathcal{A}$$

where $(\uparrow)$ is the aforementioned coercion, and $\psi$ is given by precomposition with the inclusion map from $\sigma_R(a)$ to $\sigma_{n,R}(a)$. The minor challenge comes in proving that the resulting composition so-formed is a closed embedding (given that $\phi$ itself is too). Indeed, while $(\uparrow)$ is always a closed embeddings, $\psi$ in general is not if $a$ is not invertible. Yet, we can still show that $\psi \circ (\uparrow)$ is always a closed embedding (see the proof of `isClosedEmbedding_cfcₙHom_of_cfcHom`).

## 9. autoParam AND AUTOMATION

In Section 4, we described in detail how important it was for us to unbundle the continuous functional calculus so that it is both itself a bare function, and also that it operates on bare functions. However, there is one significant downside of unbundling: to do anything non trivial, we must pass around proofs that the element $a \in \mathcal{A}$ satisfies the predicate $p$, the function $f$ is continuous on $\sigma_{\mathbb{C}}(a)$ and, in the non-unital case, $f(0) = 0$. This is cumbersome and a bit annoying, especially since oftentimes the proofs are available in the context, or easily derived from the current context, or unavailable but easily constructed.

Our solution to this conundrum is to use the `autoParam` feature of Lean. This allows the user, when writing a theorem, to specify a default tactic to use to attempt to generate a proof of one of the hypotheses. This tactic is tried when the user does not provide a proof for that hypothesis; if the tactic fails, the user must then supply a proof manually. As an example, consider the theorem:

*LeanCFC/Snippets/Automation.lean*

```
7  theorem cfcₙ_map_quasispectrum {R A : Type*} {p : A → Prop} [CommSemiring R]
8      [Nontrivial R] [StarRing R] [MetricSpace R] [IsTopologicalSemiring R] [ContinuousStar R]
9      [NonUnitalRing A] [StarRing A] [TopologicalSpace A] [Module R A] [IsScalarTower R A A]
10     [SMulCommClass R A A] [instCFCₙ : NonUnitalContinuousFunctionalCalculus R A p]
11     (f : R → R) (a : A) (hf : ContinuousOn f (σₙ R a) := by cfc_cont_tac)
12     (hf0 : f 0 = 0 := by cfc_zero_tac) (ha : p a := by cfc_tac) :
13     σₙ R (cfcₙ f a) = f '' σₙ R a := by
14   simp [cfcₙ_apply f a, cfcₙHom_map_quasispectrum (p := p)]
```

The arguments `hf`, `hf0` and `ha` each have an `autoParam` (indicated by the syntax `:= by` within the arguments), and if, when this theorem is called, they are not provided, then the tactics `cfc_cont_tac`, `cfc_zero_tac` and `cfc_tac` are used to try to construct them. Currently, each of these tactics is just a wrapper around a small collection of tactics, but in the future, we can opt to make these more sophisticated, if it is useful. As an example, this can be used to show:

*LeanCFC/Snippets/Automation.lean*

```
18 example {A : Type*} [NonUnitalCStarAlgebra A] [PartialOrder A] [StarOrderedRing A]
19     {b : A} (hb : IsSelfAdjoint b) :
20     σₙ ℝ≥0 (cfcₙ sqrt (b * b)) = sqrt '' σₙ ℝ≥0 (b * b) :=
21   cfcₙ_map_quasispectrum _ _
```

The underscores in the call to `cfcₙ_map_quasispectrum` are for the arguments to `f := sqrt` and `a := b * b`, which Lean can infer via unification in this case. The argument `ha` is constructed by the tactic `cfc_tac`, which is mostly a wrapper around the general purpose tactic `aesop`. In this case, `IsSelfAdjoint.mul_self_nonneg` is the relevant lemma marked with the `aesop` attribute, which `aesop` can then combine with the hypothesis `hb` to show that `b * b` is nonnegative (this is necessary because we're using the functional calculus over `ℝ≥0` in this example). Likewise, the argument `hf` is constructed by the tactic `cfc_cont_tac`, which is a wrapper around `fun_prop` — a general purpose tactic for proving goals related to properties of functions (e.g., continuity, differentiability, measurability, etc.). Finally, the argument `hf0` is constructed by the tactic `cfc_zero_tac`, which is also a wrapper around `aesop`; in this case it applies the `simp` lemma `NNReal.sqrt_zero`.

The net effect of this setup is that the user is able to freely apply lemmas related to the continuous function calculus, without providing proofs of the hypotheses, so long as the necessary proofs are sufficiently simple. We remark that unlike most other places in the library, the arguments `f` and `a` are intentionally left explicit, rather than implicit, despite the fact that they can be inferred from `hf`, `hf0` and `ha`. This is exactly because the latter are `autoParam`s, and since these are not provided by the user but rather autogenerated, sometimes Lean does not have enough information to infer `f` and `a`.

## 10. The isometric variation

Although we took great pains to ensure that the continuous functional calculus would be usable in contexts where a `CStarAlgebra` instance is not available, it is nevertheless the case that within the context of $C^*$-algebras, the fact that the continuous functional calculus is isometric is valuable and important. There are two ways this can be addressed. The first is to simply prove a theorem of the form `‖cfc f a‖ = ⨆ x : spectrum R a, ‖f x‖` under suitable conditions on `f` and `a : A`, in the presence of a `CStarAlgebra A` instance. The second is to provide a separate class for the isometric continuous functional calculus, which simply extends the class specified in Listing 16 by requiring that the homomorphism therein is isometric.

The latter option may seem pointless, as it will only be applicable when `A` has a metric structure, which effectively means it only applies when there is already a `CStarAlgebra` instance. However, we have indeed opted for the latter option, and there were two considerations which led us to this decision. The first is that, because there is not one (unital) continuous functional calculus, but rather three (over $\mathbb{C}$, $\mathbb{R}$ and $\mathbb{R}_{\geq 0}$), we would have to prove every theorem in triplicate; by using classes, we can have lemmas that apply to all three scalar rings simultaneously[14]. The second consideration is that, in the future, we may develop the theory of real $C^*$-algebras, and in that case, having a separate class that can apply in both the complex and real settings will be vital.

In this way, we can obtain theorems that apply to a continuous functional calculus over either $\mathbb{R}$ or $\mathbb{C}$. For instance, for every element $x$ in the spectrum of $a$, the norm of $f$ at $x$ is bounded by $\|f(a)\|$:

*LeanCFC/Snippets/Isometric.lean*

```
5 theorem norm_apply_le_norm_cfc {𝕜 A : Type*} {p : A → Prop} [RCLike 𝕜] [NormedRing A]
6     [StarRing A] [NormedAlgebra 𝕜 A] [IsometricContinuousFunctionalCalculus 𝕜 A p]
7     (f : 𝕜 → 𝕜) (a : A) ⦃x : 𝕜⦄ (hx : x ∈ spectrum 𝕜 a)
8     (hf : ContinuousOn f (spectrum 𝕜 a) := by cfc_cont_tac) (ha : p a := by cfc_tac) :
9     ‖f x‖ ≤ ‖cfc f a‖ :=
```

Because the continuous $\mathbb{R}_{\geq 0}$-valued functions don't have a norm structure in Mathlib, we sometimes state the corresponding theorems for the continuous functional calculus over $\mathbb{R}_{\geq 0}$ separately. The example below is the $\mathbb{R}_{\geq 0}$ version of the previous theorem.

*LeanCFC/Snippets/Isometric.lean*

```
12 theorem apply_le_nnnorm_cfc_nnreal {A : Type*} [NormedRing A] [StarRing A]
13     [NormedAlgebra ℝ A] [PartialOrder A] [StarOrderedRing A]
14     [IsometricContinuousFunctionalCalculus ℝ A IsSelfAdjoint] [NonnegSpectrumClass ℝ A]
15     (f : NNReal → NNReal) (a : A) ⦃x : NNReal⦄ (hx : x ∈ spectrum NNReal a)
16     (hf : ContinuousOn f (spectrum NNReal a) := by cfc_cont_tac) (ha : 0 ≤ a := by cfc_tac) :
17     f x ≤ ‖cfc f a‖₊ :=
```

## 11. Limitations and pain points

Overall, we feel that our implementation of the continuous functional calculus is quite successful. Nevertheless, there are still some rough edges that we would like to polish.

11.1. **Lack of `simp` lemmas.** Because of our choice to use bare functions everywhere, almost every single lemma (with the exception of `cfc_zero` and `cfcₙ_zero`) has some hypotheses that must be satisfied which cannot be filled by unification. As a result, these lemmas are not suitable `simp` lemmas.

For instance, the lemma `cfc_id` states that `cfc id a = a`, but it requires a proof that `a` satisfies the predicate `p` pertaining to the scalar ring `R` appearing in `id : R → R`. As outlined in Section 9, this argument is equipped with an `autoParam` which can be used to automatically construct this proof in many circumstances. One would hope that this would be sufficient to make `cfc_id` a `simp` lemma which is only applied in contexts where the proof construction succeeds. Potentially, this would be placed in a separate `simp` set, so that these lemmas (with the associated proof search induced by the `autoParams`) are not tried on every `simp` invocation that matches the pattern to prevent performance degradation.

Unfortunately, due to technical limitations (see Lean 4 issue #3475), this is not the case. That is, Lean needs the user to provide the explicit arguments for `f` and `a` in order to force Lean to generate the proofs via `autoParam`. This makes some proofs more tedious than is preferable, as one must generally provide all the rewrite steps explicitly, even the "obvious" ones. As an example, consider the following code block. Note that Lean is able to infer the explicit arguments upon using `rw`, but the underscores cannot be removed or else the `autoParam` arguments which depend on these explicit arguments will not fire and the rewrite will fail. Moreover, these cannot be used in `simp` unless the underscores are filled in, hence we must effectively rewrite one lemma at a time because, even in `simp`, we need to call `cfc_comp` twice with different explicit arguments.

[14]Actually, this is not *quite* true, as $\mathbb{R}_{\geq 0}$ poses some unique challenges and somewhat often we need separate theorems for this scalar ring. Nevertheless, this approach does significantly reduce duplication.

*LeanCFC/Snippets/Automation.lean*

```
31 example {A : Type*} [CStarAlgebra A] {a : A} [IsStarNormal a]
32     (f g : ℂ → ℂ) (hf : Continuous f) (hf : Continuous g) :
33     cfc f (cfc g a) + cfc ((2 : ℂ) * ·) (cfc g a) * cfc (f ∘ g) a =
34       cfc (fun x ↦ f (g x) + 2 * (g x) * f (g x)) a := by
35   rw [← cfc_comp _ _ _, ← cfc_comp _ _ _, ← cfc_mul _ _ _, ← cfc_add _ _ _]
36   congr! 1
```

11.2. **The headache of $\mathbb{R}_{\geq 0}$.** The use of $\mathbb{R}_{\geq 0}$ as a scalar ring is a double-edged sword. While it makes manipulation of nonnegative elements quite natural, it also introduces some complications. We often have to prove specialized versions of theorems for $\mathbb{R}_{\geq 0}$, and generic theorems don't apply. One example was given in Section 10, but we'll list a few more here. The lemma `cfc_le_iff` (and its non-unital counterpart `cfcₙ_le_iff`) shown below is only valid when the type `R` of scalars is a ring, not just a semiring. Consequently, we have to write a separate version (`cfc_nnreal_le_iff`) of this lemma for $\mathbb{R}_{\geq 0}$.

*LeanCFC/Snippets/Limitations.lean*

```
 5 theorem cfc_le_iff {R A : Type*} {p : A → Prop} [CommRing R] [PartialOrder R]
 6     [StarRing R] [MetricSpace R] [IsTopologicalRing R] [ContinuousStar R]
 7     [ContinuousSqrt R] [StarOrderedRing R] [TopologicalSpace A]
 8     [Ring A] [StarRing A] [PartialOrder A] [StarOrderedRing A] [Algebra R A]
 9     [instCFC : ContinuousFunctionalCalculus R A p] [NonnegSpectrumClass R A]
10     (f g : R → R) (a : A)
11     (hf : ContinuousOn f (spectrum R a) := by cfc_cont_tac)
12     (hg : ContinuousOn g (spectrum R a) := by cfc_cont_tac)
13     (ha : p a := by cfc_tac) :
14     cfc f a ≤ cfc g a ↔ ∀ x ∈ spectrum R a, f x ≤ g x :=
```

Likewise, the lemma `cfc_sub` (which simply states that `cfc (f - g) a = cfc f a - cfc g a`) also requires the scalars to be an actual ring. The problem is that $\mathbb{R}_{\geq 0}$ utilizes *truncated subtraction*, wherein $x - y := 0$ when $x \leq y$.

Another annoyance is related to the unitization, especially as regards the function spaces $C(s, R)$ and $C(s, R)_0$. When $0 \in s$, and $R$ is a topological *ring* (not just a topological semiring), then there is a ring isomorphism $\Phi : C(s, R) \to C(s, R)_0^{+1}$, where the unitization is taken over $R$. Here, $\Phi(f) := (f(0), f - f(0))$, and $\Phi^{-1}(a, g) := a + g$. However, this isomorphism clearly falters when $R$ is only a semiring, these rings are not isomorphic in that case.

## 12. Simplified statements and proofs

Up to this point, we have been expounding upon our design, and explaining why our version consistutes a better definition than the alternatives. In this section, we will provide some brief examples comparing statements and proofs using our interface with those using a different one. We've tried to make these comparisons *fair* in the sense that, although we have developed a large API for our preferred interface, we limit ourselves to only using that which corresponds to API for the mocked-up alternatives.

12.1. **Monomorphism as the primary interface.** Our first comparison is to a version of the continuous functional calculus where the $*$-monomorphism `C(spectrum ℂ a, ℂ) →⋆ₐ[ℂ] A` is the primary interface. One can think of this as being somehow halfway in-between Listing 1 and Listing 2 in that we remove the dependent typing on the codomain, but leave in place the bundled structure of the continuous maps and the $*$-homomorphism. In that sense, this should be almost uniformly more pleasant to work with than Listing 1 without losing many of the benefits, but we will show it still leaves us wanting in comparison to Listing 2.

Consider a statement as simple as $f(-a) = (x \mapsto f(-x))(a)$ for $f$ continuous on the spectrum of $-a$, which was one of our litmus tests. In the monomorphism approach, this is stated and proved (with some limited background API) as:

*LeanCFC/Snippets/Monomorphism.lean*

```
109 theorem cfcMono_map_neg {a : A}
110     [IsStarNormal a] (f : C(spectrum ℂ (-a), ℂ)) :
111     cfcMono (-a) f =
112       cfcMono a (f.comp
113         ⟨Subtype.coind (- · : spectrum ℂ a → ℂ) (by simp), by fun_prop⟩) := by
114   refine Eq.symm <| cfcMono_map_comp_of_eq (a := a) (b := -a) (g := f)
115     (f := ⟨(- ·), by fun_prop⟩) ?_ (fun _ ↦ rfl)
116   conv_lhs => rw [← cfcMono_map_id a, ← map_neg]
117   rfl
```

Note the complexity of the statement due to the need to bundle the map $x \mapsto -x$ into a continuous function from $\sigma_{\mathbb{C}}(a)$ to $\sigma_{\mathbb{C}}(-a)$. Whereas in our approach, the statement matches what one would write on paper, and the proof is a bit simpler. What follows is the same statement as Listing 11, but we note that the proof is now even simpler because we make use of the `autoParam` arguments as described in Section 9.

*LeanCFC/Snippets/Monomorphism.lean*

```
135 lemma cfc_map_neg {a : A} [IsStarNormal a] (f : ℂ → ℂ)
136     (hf : ContinuousOn f (spectrum ℂ (-a))) :
137     cfc f (-a) = cfc (fun x ↦ f (-x)) a := by
138   rw [← spectrum.neg_eq] at hf
139   rw [cfc_comp' f (- ·) a (hg := by simpa), cfc_neg, cfc_id' ..]
```

However, when one ventures just slightly beyond the most straightforward statements, the complexity of the monomorphism approach explodes, whereas our approach remains manageable. Consider, for example, the statement $(x \mapsto f(-x))((-a)^3) = f(a^3)$ for $f$ continuous on the spectrum of $a^3$. In the monomorphism approach, the statement itself is again marred by the need to bundle the map $x \mapsto -x$ into a continuous function from $\sigma_{\mathbb{C}}(a^3)$ to $\sigma_{\mathbb{C}}((-a)^3)$; the proof is brittle (e.g., we were unable to easily rewrite with the lemma `cfcMono_congr_apply` and were forced to use `convert` instead); and the proof is simply longer than would be preferred.

*LeanCFC/Snippets/Monomorphism.lean*

```
260 example {a : A} [IsStarNormal a] [Nontrivial A] (f : C(spectrum ℂ (a ^ 3), ℂ)) :
261     cfcMono ((-a) ^ 3) (f.comp
262       ⟨Subtype.coind (- · : spectrum ℂ ((-a) ^ 3) → ℂ)
263         (by simp [Odd.neg_pow (n := 3) (by decide), ← spectrum.neg_eq]), by fun_prop⟩) =
264     cfcMono (a ^ 3) f := by
265   rw [← cfcMono_map_neg' f]
266   convert (cfcMono_congr_apply ?_ _).symm using 2
267   · ext
268     rfl
269   · grind
```

Compare this with our approach, where the statement is clear, and the proof consists of only two simple rewrites. Moreover, our automation is capable of filling in the necessary proofs of continuity and normality behind the scenes so that the user need not even worry about it.

*LeanCFC/Snippets/Monomorphism.lean*

```
273 example {a : A} [IsStarNormal a] (f : ℂ → ℂ) (hf : ContinuousOn f (spectrum ℂ (a ^ 3))) :
274     cfc (fun x ↦ f (-x)) ((-a) ^ 3) = cfc f (a ^ 3) := by
275   rwa [Odd.neg_pow (by decide), cfc_map_neg']
```

To truly drive home the point, let's up the complexity of the statement ever so slightly more. In this case, we would like to show that $(x \mapsto f(\overline{x}))((-a)^3) = (x \mapsto f(-x))((a^*)^3)$, for $f$ continuous on $\sigma_{\mathbb{C}}((-a^*)^3)$. Here the monomorphism approach balloons. The statement requires extracting some proof terms into `haveI` statements in order to be even moderately readable. The proof is brittle and long: we effectively make the goal state more and more complicated as we rewrite until we finally manage to get to a place where we can prove two functions are equal. We encourage the reader to peruse the goal state immediately preceding the line `congr! 1` in the next example for context.

*LeanCFC/Snippets/Monomorphism.lean*

```
291 example {a : A} [IsStarNormal a] (f : C(spectrum ℂ ((-star a) ^ 3), ℂ)) :
292     haveI h₁ (x : spectrum ℂ ((-a) ^ 3)) : star ↑x ∈ spectrum ℂ ((-star a) ^ 3) := by
293       simp only [← star_neg, ← star_pow, spectrum.map_star, Set.star_mem_star]
294       exact x.2
295     haveI h₂ (x : spectrum ℂ (star a ^ 3)) : -↑x ∈ spectrum ℂ ((-star a) ^ 3) := by
296       simp only [Odd.neg_pow (n := 3) (by decide), ← spectrum.neg_eq, Set.neg_mem_neg]
297       exact x.2
298     cfcMono ((-a) ^ 3) (f.comp
299       ⟨Subtype.coind (star · : spectrum ℂ ((-a) ^ 3) → ℂ) h₁, by fun_prop⟩) =
300     cfcMono (star a ^ 3) (f.comp
301       ⟨Subtype.coind (- · : spectrum ℂ (star a ^ 3) → ℂ) h₂, by fun_prop⟩) := by
302   conv_rhs => rw [← cfcMono_map_neg']
303   conv_lhs => rw [← cfcMono_map_star']
304   rw [cfcMono_congr_apply (b := star ((-a) ^ 3)) (a := -star a ^ 3) ?h_eq]
305   case h_eq => simp [Odd.neg_pow (n := 3) (by decide)]
306   congr! 1
307   ext
308   simp [Subtype.coind, Homeomorph.setCongr]
```

In contrast, in our approach, the statement is again clear. Note how the proofs which were required for the *statement* in the previous example are instead deferred into the *proof* in the following example; this is the "write first, think later" principle in action that we discussed in Section 4. The proof is slightly more involved than in previous problems, as this one has more moving parts, but it is still manageable.

*LeanCFC/Snippets/Monomorphism.lean*

```
314 example {a : A} [IsStarNormal a] (f : ℂ → ℂ)
315     (hf : ContinuousOn f (spectrum ℂ ((-star a) ^ 3))) :
316     cfc (fun x ↦ f (star x)) ((-a) ^ 3) = cfc (fun x ↦ f (-x)) (star a ^ 3) := by
317   rw [cfc_comp_star (hf := ?hf₁), cfc_comp_neg (hf := ?hf₂)]
318   case hf₁ => rwa [Set.image_star, ← spectrum.map_star, star_pow, star_neg]
319   case hf₂ => rwa [Set.image_neg_eq_neg, spectrum.neg_eq, ← Odd.neg_pow (by decide)]
320   simp [Odd.neg_pow (n := 3) (by decide)]
```

If, in the preceding example problem, one considers instead a globally continuous function $f$ (i.e., continuous on all of $\mathbb{C}$), then the statement in the monomorphism approach is correspondingly simpler, but in our approach nothing needs to be changed except the hypothesis. Moreover, the proof in the monomorphism approach is exactly the same; there is no simplification. Whereas in our approach, for globally continuous $f$ the automation can more effectively reduce our workload, thereby simplifying the proof. See the relevant examples in *LeanCFC/Snippets/Monomorphism.lean*.

12.2. **Allowing only complex-valued functions.** We emphasized the importance of, and went to some lengths to incorporate, $\mathbb{R}$ (and even $\mathbb{R}_{\geq 0}$) as scalar rings for the continuous functional calculus, so that are functions may be defined on and take values there instead. In this example, we will only use the interface we designed for the continuous functional calculus, having established in the previous section its superiority. In the first example we will tie our hands behind our backs by only allowing complex-valued functions and the continuous functional calculus over $\mathbb{C}$. The theorem proven here is one that actually occurs in practice (cf. `SpectrumRestricts.nnreal_iff_nnnorm`), but we have slightly modified the statement to better exhibit the differences between working only over $\mathbb{C}$, or allowing $\mathbb{R}$ as well.

*LeanCFC/Snippets/SpectrumRestricts.lean*

```
30 example [Nontrivial A] {a : A} (ha : IsSelfAdjoint a) {t : ℝ} (ht : ‖a‖ ≤ t)
31     (h : ∀ x ∈ spectrum ℂ a, 0 ≤ x) : ‖algebraMap ℝ A t - a‖ ≤ t := by
32   rw [IsScalarTower.algebraMap_apply ℝ ℂ A, Complex.coe_algebraMap]
33   rw [← cfc_id' ℂ a, ← cfc_const (t : ℂ) a, ← cfc_sub ..]
34   refine IsGreatest.norm_cfc (t - · : ℂ → ℂ) a |>.isLUB.2 ?_
35   rintro - ⟨x, hx, rfl⟩
36   simp only
37   lift x to ℝ using (h x hx).2.symm
38   rw [← Complex.ofReal_sub, Complex.norm_real, Real.norm_eq_abs]
39   have hx' : 0 ≤ x := by simpa using h x hx
40   have hxt : x ≤ t := by
41     simpa [abs_of_nonneg hx'] using (spectrum.norm_le_norm_of_mem hx |>.trans ht)
42   grind
```

For comparison, here is nearly the same statement and proof, but using the continuous functional calculus over $\mathbb{R}$ as well. Notice how there is simply less friction and clutter in approaching it this way. We realize that one hypothesis is slightly different, but in the code artifact we provide a (short!) proof that the differing hypotheses are equivalent in the presence of the other hypotheses.

*LeanCFC/Snippets/SpectrumRestricts.lean*

```
44 example [Nontrivial A] {a : A} (ha : IsSelfAdjoint a) {t : ℝ} (ht : ‖a‖ ≤ t)
45     (h : ∀ x ∈ spectrum ℝ a, 0 ≤ x) : ‖algebraMap ℝ A t - a‖ ≤ t := by
46   rw [← cfc_id' ℝ a, ← cfc_const t a, ← cfc_sub ..]
47   refine IsGreatest.norm_cfc (t - ·) a |>.isLUB.2 ?_
48   rintro - ⟨x, hx, rfl⟩
49   have := spectrum.norm_le_norm_of_mem hx |>.trans ht
50   grind [Real.norm_eq_abs, h x hx]
```

While the difference in difficulty is not as stark as with the examples in Section 12.1, which is primarily due to the fact that the *statements* are now manageable in both cases, the kind of friction witnessed here can accumulate as one develops more theory.

## 13. Applications

Following the addition of our interface to Mathlib, multiple contributions have relied on it to further develop the mathematical aspects of the theory. We interpret this as an indication that, despite the limitations outlined in Section 11, our interface is convenient to work with even in cases where the computations involved are quite convoluted. While a detailed discussion of any of these contributions is outside the scope of this paper, we list and comment on some of them here.

- **Existence of approximate units in non-unital $C^*$-algebras (Jireh Loreaux, 2024).** This relies on quite heavy calculations using the continuous functional calculus.
  See lemma `CStarAlgebra.increasingApproximateUnit`.
- **Continuity of the functional calculus in both variables (Jireh Loreaux, 2025).** Let $\mathcal{A}$ be a non-unital $C^*$-algebra and $S$ a compact subset of $\mathbb{C}$. Then, the map $(f, a) \mapsto f(a)$ is continuous on the set $X$ of all pairs $(f, a)$ such that $a \in \mathcal{A}$ is normal, $\sigma_{\mathbb{C}}(a) \subseteq S$ and $f : \mathbb{C} \to \mathbb{C}$ is continuous on $\sigma_{\mathbb{C}}(a)$. Here, the topology on $X$ is that of uniform convergence on $S$ of the first component, and norm-convergence in $\mathcal{A}$ of the second component. We note that the ability to write `cfc f a` quite liberally contributes quite significantly to the readability and simplicity of this statement in Lean.
  See lemmas `continuousOn_cfc_setProd` and `continuousOn_cfcₙ_setProd` for example.
- **The maps $a \mapsto |a|$, $a \mapsto a^p$ and $a \mapsto \log(a)$ on elements of a non-unital $C^*$ algebra are continuous on their respective domains (Jon Bannon, Frédéric Dupuis and Jireh Loreaux, 2024–2026).**
  See lemmas `CFC.continuous_abs`, `CFC.continuousOn_nnrpow` and `CFC.continuousOn_log`.

- **The maps $a \mapsto a^p$ $(0 \leq p \leq 1)$ and $a \mapsto \log(a)$ on positive elements of a non-unital $C^*$ algebra are concave and non-decreasing (Frédéric Dupuis, 2025).** See lemmas `CFC.monotone_nnrpow` and `CFC.log_monotoneOn`.
- **Projections in a non-unital $C^*$-algebra are the extreme points of the set of positive elements of norm at most $1$ (Jon Bannon, Jireh Loreaux and Monica Omar, 2026).** See lemma `isStarProjection_iff_mem_extremePoints_setOf_nonneg_inter_unitClosedBall`[15]
- **In a unital $C^*$-algebra, unitary elements span the algebra and the unitary group is locally path connected (Jireh Loreaux, 2025).** See lemmas `CStarAlgebra.span_unitary` and `Unitary.instLocPathConnectedSpace`.

## 14. Future Work

**14.1. A `cfc` manipulation tactic.** By perusing the code artifact associated to this paper, the reader may notice that the proofs often involve rewrites (or `simp` calls) that rewrite in the *reverse* direction (i.e., with `←`) of the lemma statement. From a mathematical perspective, the reason is clear: in order to show that two expressions involving the continuous functional calculus are equal, it suffices to write both as applications of the continuous functional calculus for some element $a \in \mathcal{A}$ to some functions $f, g : \mathbb{C} \to \mathbb{C}$ and then show that $f$ and $g$ agree on the spectrum of $a$. This last step is the lemma `cfc_congr` or `cfcₙ_congr` that appears repeatedly. And of course, this is a natural thing to do because working with functions is simpler than working with elements in the $C^*$-algebra.

The reader may wonder why we state all these lemmas in the wrong direction. While this is something we could switch, and maybe we should do so at some point, there is a reason it ended up this way. Throughout Mathlib, when $\phi : \mathcal{A} \to \mathcal{B}$ is a morphism of semigroups (or anything stronger, such as a $*$-homomorphism), the lemma `map_mul` applies, which states that $\phi(xy) = \phi(x)\phi(y)$ for $x, y \in \mathcal{A}$. Moreover, this lemma is a `simp` lemma, and so it is applied automatically whenever that tactic is called. Recalling that `cfc` itself (or more precisely, `cfcHom`) is a $*$-homomorphism, we can see that `cfc_mul` is just a special case of `map_mul`.

To emphasize, lemmas in the library which apply to morphisms and are marked `simp` (e.g., `map_mul`) generally push the morphism to the *leaves* of the expression (viewed as tree). In contrast, when applying the analogous lemmas (such as `cfc_mul`) for the continuous functional calculus, we generally want to pull `cfc` to the *head* of the expression, so it's exactly the reverse of what we do normally. Of course, often we *really do* want the lemmas in the direction currently stated (e.g., `cfc_id`), and we would like it if they were `simp` lemmas (but only in this direction). At the same time, when moving `cfc` to the head of an expression, even `cfc_id` is used in the opposite direction.

This suggests that what we truly need is a tactic that can manipulate expressions involving the continuous functional calculus in an intelligent way. Such a tactic should be able to pull `cfc` to the head of an expression, or push it to the leaves as necessary. As it applies lemmas, it should collect goals that cannot be solved by the existing automation, and leave them for the user to prove. This would reduce the burden on the user to continually provide the entire list of rewrites, and, if written properly, would avoid the need for many targeted rewrites like the invocation `nth_rw 2 [← cfcₙ_id' ℝ a]`.

**14.2. Real $C^*$-algebras.** Throughout the development process, we have been careful to ensure that our implementation is as general as possible. In particular, we want this to be usable for real $C^*$-algebras when those eventually enter Mathlib. Note that unlike most areas of mathematics where the real theory may generally be developed in parallel, or even prior to, the complex theory, in the context of $C^*$-algebras, the complex theory must be developed first. Indeed, even the spectrum of an element in a real $C^*$-algebra is a subset of the complex plane (and not just a subset of the real line). One way to understand the reason for this is that the category of real $C^*$-algebras is equivalent to the category of complex $C^*$-algebras equipped with a conjugate-linear multiplicative involution (or equivalently, a linear antimultiplicative involution).

Our implementation of the continuous functional calculus over $\mathbb{R}$ and $\mathbb{R}_{\geq 0}$ should work out-of-the-box for real $C^*$-algebras, once the relevant instances are supplied. On the other hand, the continuous functional calculus over $\mathbb{C}$ (for normal elements) is markedly different in the case of real $C^*$-algebras. In particular, it is *not* a $*$-homomorphism from $C(\sigma_{\mathbb{C}}(a), \mathbb{C})$ into $\mathcal{A}$, but the domain is rather the subalgebra of those continuous functions which commute with complex conjugation (i.e., $f(\overline{z}) = \overline{f(z)}$). For example, if $a \in \mathcal{A}$ is a normal element in a unital real $C^*$-algebra, then even though the function $z \in \mathbb{C} \mapsto e^{iz}$ is continuous everywhere (hence on the spectrum of $a$), in general $e^{i\overline{z}} \neq \overline{e^{iz}}$, and hence it is not possible to make sense of $e^{ia}$ via the continuous functional calculus for real $C^*$-algebras. Given that this is sufficiently different from all the other cases, our current attitude is that this should be handled via a bespoke interface, even if the proofs ultimately boil down to an appropriate application of the continuous functional calculus over $\mathbb{C}$ on a certain complex $C^*$-algebra.

[15]This result is located on a newer version of Mathlib than anything else listed in this paper, which we have kept consistent at v4.28.0.

## Code artifact

Code associated with this paper is available at `https://github.com/j-loreaux/LeanCFC/releases/tag/v2`

## Acknowledgements

The authors thank Frédéric Dupuis for many thoughtful discussions during the development of this project, and for much code review upon submission to Mathlib. More generally, we're indebted to the entire Mathlib and Lean communities, without whose efforts this work would not have been possible.

In addition, we would like to thank one of the anonymous referees for their careful reading and extensive feedback on an earlier version of this manuscript.

Université Paris Cité, Sorbonne Université, CNRS, IMJ-PRG, F-75013 Paris, France

Southern Illinois University Edwardsville, Edwardsville, Illinois 62026, USA